\documentclass[a4paper,twoside,10pt]{article}
\usepackage{amsmath,amssymb,amsfonts} 
\usepackage{stackrel}
\usepackage{amssymb}
\usepackage{graphics}                 
\usepackage{enumerate}
\usepackage{amsthm}
\usepackage{graphicx}
\usepackage{caption}
\usepackage{subcaption}
\usepackage{hyperref}
\usepackage{multirow}
\usepackage{eufrak}
\usepackage[numbers,sort&compress]{natbib}
\usepackage{marginnote}
\usepackage{color}

 \newcommand {\R}{\mathbb R}
 \newcommand {\T}{\mathbb T}
 \newcommand {\sS}{\mathbb S}
\newcommand {\C}{\mathbb C}
 \newcommand {\Z}{\mathbb Z}
 \newcommand {\N}{\mathbb N}
 \newcommand {\e}{\mathrm e}
 \renewcommand {\i}{\mathrm i}
\newcommand {\Fix}{\mathrm{Fix}}
 \renewcommand {\S}{\mathbf{S}}
 \renewcommand {\d}{\mathrm{d}}
 \newcommand {\D}{\mathrm{D}}
\newcommand {\DD}{\mathbf{D}}
 \newcommand {\kC}{\mathcal{C}}
\newcommand {\kX}{\mathcal{X}}
\newcommand {\Mat}{\text{Mat}}

\graphicspath{{figures/}}
\allowdisplaybreaks[2]
\date{\today}
\title{Symmetric bifurcation analysis of synchronous states of
 time-delayed coupled Phase-Locked Loop oscillators}
\author{ \textbf{Diego Paolo Ferruzzo Correa}\\
Universidade de S\~ao Paulo, Escola Polit\'ecnica,\\
Departamento de Telecomunica\c c\~oes e Controle\\
S\~ao Paulo-SP, Brasil\\
dferruzzo@usp.br\\{}\\
\textbf{Claudia Wulff}\\
Department of Mathematics, University of Surrey, UK\\
c.wulff@surrey.ac.uk\\{}\\
\textbf{Jos\'e Roberto Castilho Piqueira}\\
Universidade de S\~ao Paulo, Escola Polit\'ecnica, \\
Departamento de Telecomunica\c c\~oes e Controle\\
S\~ao Paulo-SP, Brasil\\
piqueira@lac.usp.br}
\theoremstyle{definition}

\newtheorem{lemma}{Lemma}[section]

\newtheorem{definition}{Definition}[section]

\theoremstyle{remark}
\newtheorem*{remark}{Remark}

\numberwithin{equation}{section}
\newcommand{\Arctan}{\mathrm{Arctan}}
\renewcommand{\Re}{\mathrm{Re}}
\renewcommand{\Im}{\mathrm{Im}}

\newcommand{\note}[1]{\vspace{5 mm}\par \noindent
  \marginpar{\textsc{Note}} \framebox{\begin{minipage}[c]{0.95
        \textwidth} \flushleft \tt #1 \end{minipage}}\vspace{5
    mm}\par}
\renewcommand{\note}[1]{}
\renewcommand{\marginnote}[1]{}
\renewcommand{\fbox}[1]{}
\begin{document}
\maketitle
\begin{abstract}
In recent years there has been an increasing interest in studying
time-delayed coupled networks of oscillators since these occur in many
real life applications. In many cases symmetry patterns can emerge in
these networks, as a consequence a part of the system might repeat
itself, and properties of this subsystem are representative of the
dynamics on the whole phase space. In this paper an analysis of the
second order N-node time-delay fully connected network is presented which is based
on previous work by Correa and Piqueira \cite{Correa2013} for a 2-node network. This study is carried out using symmetry groups. We show the existence of multiple eigenvalues forced by symmetry, as well as the existence of Hopf bifurcations.
 Three different models are used to analyze the network dynamics,
 namely, the full-phase, the phase, and the phase-difference model. We
 determine a finite set of frequencies $\omega$, that might correspond
 to Hopf bifurcations in each case  for critical values of the
 delay. The $S_n$ map  is used to actually find Hopf bifurcations along with numerical calculations using the Lambert W function. Numerical simulations are used in order to confirm the analytical results.
Although we restrict attention to second order nodes, the results could be  extended to higher order networks provided the time-delay in the connections between nodes remains equal.
\end{abstract}

\tableofcontents


\section{Introduction}
\label{sec:intro}
Coupled oscillators present a great variety of interesting phenomena and provide models for many different areas in engineering, biology, chemistry, economy, etc. There is a considerable body of literature on coupled oscillators, see e.g.~\cite{Bueno2010, Carareto2009986, Correa2013, Piqueira2006, Piqueira2011, Piqueira2003}, in particular on different network configurations without time-delay. In~\cite{Alexander1986} a global bifurcation analysis for a network of linear coupled oscillators without delay was presented with applications in chemical processes, and in~\cite{Wu1998} an extension of this work was presented considering the lag among nodes as bifurcation parameter for neural networks with symmetry; an analysis of several configurations of oscillators  with smooth coupling functions is presented in~\cite{Jia2013, Kazanovich2013}; similar results considering patterns emerging in networks of coupled oscillators with time-delay can be found in~\cite{Yao2013}. 
         
We are interested in obtaining the simplest model for an N-oscillator
 second-order network with time-delay between oscillators; for
 this purpose we shall choose a Phase-Locked Loop (PLL) as node,
 see~\cite{Best2007}; the main difference between a PLL and other
 kinds of oscillators used frequently in the literature is that a
 PLL can oscillate by itself and its frequency can be controlled by an external signal, see~\cite{Floyd2005}. In order to obtain a proper mathematical model for a single node we shall take as a starting point the classical approach as presented by Floyd~\cite{Floyd2005} and Kudrewicz~\cite{JacekKudrewicz2007}, and for the network we will use the model introduced by Piqueira-Monteiro~\cite{Piqueira2006}, but here additionally we shall compare three different models, namely, the full-phase model, the phase model and the phase-difference model.
 
Numerical results obtained for these models are used in order to validate our analytical conclusions, especially when discussing bifurcation points, which is the main aim of our research.

The structure of this paper is as follows: 
In section~\ref{sec:Full_phase_model} the full-phase model for an N-node fully connected time-delay network is reviewed, the focus is the symmetry of the network and to find irreducible representations, bifurcations are analyzed in each  of the isotypic components. In sections~\ref{sec:Phase_model} and~\ref{sec:phase_diff_model} a comparative analysis between the phase model and the phase-difference model is performed using the results obtained for the full-phase model, and finally  conclusions and insights for future research are presented in section~\ref{sec:conclusions}.



\section{Full-phase model}
\label{sec:Full_phase_model}
In~\cite{Piqueira2006} a model for a fully connected N-node network
with time-delay is presented, each node is a second-order PLL
oscillator, see~\cite{Best2007, Floyd2005}; in this model the so called "double-frequency" term is neglected arguing that its influence is  suppressed by the local dynamics in each node, see~\cite{Acebron2005, Bueno2010, Correa2013, Floyd2005, JacekKudrewicz2007, Monteiro2003, Piqueira2006, Piqueira2003}. In what follows we will use this model but include the double-frequency term, thus we have that equation $(3.9)$ in~\cite{Piqueira2006} becomes 
\begin{equation}
  \label{e.fullPhase}
\begin{array}{l}
  \ddot{\phi}_i(t)+\mu\dot{\phi}_i(t) -\mu\omega_M = \\[1ex]
\quad\quad \dfrac{K\mu}{N-1}\displaystyle\sum_{\stackrel[j\neq i]{j=1}{}}^N\left[\sin(\phi_j(t-\tau)-\phi_i(t))+\sin(\phi_j(t-\tau)+\phi_i(t))\right],
\end{array}
\end{equation}
$i=1,\ldots, N$. The ``double-frequency'' term is embedded in the term $\sin(\phi_j(t-\tau)+\phi_i(t))$. This equation models the dynamics for the $i$-th oscillator in the N-node network; we call $\phi_i(t)$ 
\begin{equation*}
  \phi_i(t):=\theta_i(t)+\omega_Mt,~~~i=1,\ldots,N,
\end{equation*}
the full-phase  of the $i$-th oscillator where $\theta_i(t)$ is the
local instantaneous phase, $\omega_M>0$ represents the local frequency, $K>0$ and $\mu>0$ are called gains, and $\tau\geq 0$ is the time-delay.

Note that \eqref{e.fullPhase} has  equilibria at $\phi_i = \phi^{\pm}$, where
\begin{equation}
\label{eq:full_phase_equilibria_points}
\begin{array}{l}
2\phi^+=\arcsin\left(-\dfrac{\omega_M}{K}\right)~{\rm~mod}~ 2\pi,\\\\
 2\phi^-=\pi-\arcsin\left(-\dfrac{\omega_M}{K}\right)~{\rm~mod ~2\pi}
\end{array}
\end{equation}
and $\arcsin(\cdot)$ takes its values in $[-\pi/2,\pi/2]$.
Here we restrict to $\omega_M/K\leq1$
to ensure existence of equilibria.
Note that when $\tilde K=\omega_M/K=1$ then $\phi^+=\phi^-$
if we choose $\phi^+,\phi^- \in (-\pi,\pi]$; moreover the curve $\phi(\tilde K)$  of $\S_N$-invariant equilibria 
parametrized by $\tilde K$ with $\phi(\tilde K) = \phi^\pm (\tilde K)$ has a  saddle node bifurcation at $\tilde K=1$.


\subsection{$\S_N$-symmetry and irreducible representations in the full-phase model}
\label{subsec:SN_symm_full_phase_model}
We now show that \eqref{e.fullPhase} has $\S_N$-symmetry, where $\S_N$ is the group of 
all permutations $\gamma$ of $N$ elements.
A differential equation $\dot X(t) = F(X(t))$ posed on a phase space $\kX$ is equivariant with respect to the action of
a Lie group  $\Gamma$ on $ \kX$ if 
\begin{equation*}
\gamma F(X) = F(\gamma X) ~~\mbox{for all}~~ X \in \kX, \gamma \in \Gamma,
\end{equation*}
see \cite{Golubitsky2002, Golubitsky1988}.
In this case the phase space is 
$\kX = \kC([-\tau,0),\R^{2N})$, the Banach space of continuous functions from $[-\tau,~0]$ into $\R^{2N}$ 
equipped with the usual supremum norm 
\begin{equation}
\|x\|=\stackrel[-\tau\leq\theta\leq0]{}{\text{sup}}|x(\theta)|,~~~x\in\kC([-\tau,0],\R^{2N}),
\label{eq:norm}
\end{equation} see, e.g., \cite{Krawcewicz1999}; we write $\phi = (\phi_1,\ldots, \phi_N)$,
with $\phi_j \in \kC([-\tau,0),\R)$, $j=1,\ldots, N$,
and let  $x= (x^{(1)}, \ldots, x^{(N)}) \in \kX$ where $x^{(i)} = (x^{(i)}_1,x^{(i)}_2)$, 
and  $x^{(i)}_1=\phi_i$ and $x^{(i)}_2=\dot{\phi}_i$, $i=1,\ldots, N$.
If $x:[-\tau,A]\to\R^n$ is a continuous function with $A>0$ and if $t\in[0,A]$, then $X(t)\in\kC([-\tau,0],\R^n)$
 is defined by
\begin{equation}
X(t)(\theta)=x(t+\theta),~~~\theta\in[-\tau,0], t \in [0,A].
\label{eq:x_t_theta}
\end{equation}
Then \eqref{e.fullPhase} takes the form
\begin{equation*}
\frac{\d}{\d t} X(t) = F(X(t),\eta)
\end{equation*}
where 
\begin{equation*}
  \label{e.F}
  (F(x))(\theta)=\left\{
      \begin{array}{ccc} \frac{\partial x}{\partial \theta}(\theta)
        &,&-\tau\leq\theta\leq 0\\
        f(x(0),x(-\tau),\eta)&,&\theta=0.
      \end{array} 
\right.
\end{equation*}
Here $\eta = (\mu,K,\omega_M, \tau) \in \R^4$ is a parameter and
 $f= (f^{(1)}, \ldots, f^{(N)})$ is such that \eqref{e.fullPhase} can also be rewritten as 
  autonomous nonlinear delay differential equation (DDE)
\begin{equation}
 \dot{x} =f(x,x_\tau,\eta),
\label{eq:RFDE}
\end{equation}
i.e.,  $f^{(i)}= (f^{(i)}_1, f^{(i)}_2)$, where
\begin{equation}
\label{e.f}
\begin{array}{lcl}
  f^{(i)}_1(x, x_\tau) &=  & x_2^{(i)}, \\
        f^{(i)}_2(x, x_\tau)  &=    &              -\mu x_2^{(i)}+\mu\omega_M\\
&& +\dfrac{K\mu}{N-1}\displaystyle\sum_{\stackrel[j\neq i]{j=1}{}}^N\left[\sin(x_{1,\tau}^{(j)}-x_1^{(i)}) +\sin(x_{1,\tau}^{(j)}+x_1^{(i)})\right].
\end{array}
\end{equation}
Here, for short we write $X(t)(-\tau)=x(t-\tau)=:x_\tau$ and $\Gamma = \S_N$  acts on $\kX$  via 
\begin{equation}
(\gamma x)(\theta)= \gamma x(\theta),~~~\theta\in[-\tau,0].
\label{eq:symmetry_2}
\end{equation}
Note that $\S_N$  is generated by the transpositions $\pi_{ij}\in\Mat(2N)$, which swap $x^{(i)}$ with $x^{(j)}$.
 To show that  \eqref{e.fullPhase} has  $\S_N$-symmetry it is thus sufficient to prove that 
  $f\circ\pi_{ij}=\pi_{ij}\circ f$ for all $\pi_{ij}$. We compute that 
\begin{equation}
\begin{array}{rcl}
  (f\circ\pi_{ij})_2^{(i)}(x ,x_\tau) &= &
                                -\mu x_2^{(j)} +\mu\omega_M\\
&& +\dfrac{K\mu}{N-1}\displaystyle\sum_{\stackrel[\ell\neq j]{l=1}{}}^N\left[\sin(x_{1,\tau}^{(\ell)} -x_1^{(j)})+\sin(x_{1,\tau}^{(\ell)} +x_1^{(j)} )\right]
            \\               
                      &=&f^{(j)}_2(x ,x_\tau)\\
                      &=&(\pi_{ij}\circ f)_2^{(i)}(x ,x_\tau);       
\end{array}
\end{equation}
Since $\pi_{ij}=\pi_{ji}$, this  argument also gives $(f\circ\pi_{ji})^{(j)}=(\pi_{ji}\circ f)^{(j)}$, and, since  for all $k\neq i,j$ we have $(\pi_{ij}\circ f)^{(k)}=f^{(k)}= (f\circ\pi_{ij})^{(k)}$, we see that $f$ and $\pi_{ij}$ commute for all $i,j$
which proves $\S_N$-symmetry of  \eqref{e.fullPhase}.

A space $V$ is called $\Gamma$-invariant if $gV\subseteq V$, for all $g\in\Gamma$.
When a compact group acts on a space $V$, we can decompose the space into $\Gamma$-invariant subspaces of smaller dimension. The smallest blocks for such a decomposition are said to be irreducible.
When $\Gamma$ is finite then  there is a finite number of distinct $\Gamma$-irreducible subspaces of $V$, call these $U_1,\ldots,U_t$. Define 
$V_k$ to be the sum of all $\Gamma$-irreducible subspaces $U$  of $V$ such that $U$ is $\Gamma$-isomorphic to $U_k$. Then
\begin{equation}
V=V_1\oplus\cdots\oplus V_t.
\end{equation}
This generates a unique decomposition of $V$ into the so-called \textit{isotypic} components $V_j$ of $V$~\cite{Golubitsky1988}.
A representation of a group $\Gamma$ on a vector space $V$ is said to be absolutely irreducible if the only linear mapping on $V$ that 
commutes with all  $\gamma \in \Gamma$ is a  scalar multiple of the identity~\cite{Golubitsky1988}.

The $\S_N$-symmetry acting on 
$\R^{N}$ by permuting coordinates is called 
 the permutation representation. This has the trivial subrepresentation consisting of vectors whose coordinates are all equal. The orthogonal complement consists of those vectors whose coordinates sum to zero, and when $N\geq2$, the representation on this subspace is an $N-1$-dimensional absolutely irreducible representation of $\S_N$, called the standard representation, see e.g. \cite{Dias2009}.  
In other words, $\R^{N}$ decomposes as 
\begin{equation}\label{e.U}
\R^{N} = 
\Fix(\S_N)\oplus U,
\end{equation}
 where 
for any   subgroup $H$ of  the action of a Lie group $\Gamma$  on $\kX$ the fixed-point subspace of $H$ is given by
\begin{equation*}
\Fix_{\kX}(H)= \Fix(H) = \{ X\in\kX~\mid~hX=X,~\forall h\in H\} 
\end{equation*}%
and $U= (\Fix(\S_N))^\bot \cong\R^{N-1}$ is $\S_N$-invariant and irreducible.
Moreover 
$\R^{2N}$, the phase space for $(\phi,\dot\phi)$, decomposes as $\R^{2N}=\Fix(\S_N) \oplus V$ where 
$\Fix(\S_N)=\R^2$, and $V=U\oplus U$ are the isotypic components of the $\S_N$-action on $\R^{2N}$. 

Note that if $A$ is a linear operator on a vector space $V$ with a Lie
group $\Gamma$ acting linearly on $\kX$ and
$A$ is
 $\Gamma$-equivariant, i.e., $\gamma A=A \gamma$ for all $\gamma \in \Gamma$ then $A$ has a block decomposition, more precisely, $A(V_j)\subseteq V_j$ for all isotypic components $V_j$ of $V$. Moreover if $V_j$ is the isotypic component of
an absolutely irreducible representation $U_j$ of dimension $n$ then $A|_{V_j}$ consists of $n$ identical blocks.

The equilibria  $x^\pm$ with $(x_1^{(i)},x_2^{(i)})=(\phi^{\pm},0)$, $i=1,\ldots, N$,  from \eqref{eq:full_phase_equilibria_points}
are $\S_N$-invariant and hence their linearization $A(\eta) = \D F(x^\pm,\eta) $ with $F$ as in \eqref{e.F}
is $\S_N$-equivariant. Note that
\begin{equation*}
(A(\eta) x)(\theta)= \left\{
    \begin{array}{ccc}
      \dfrac{dx}{d\theta}&,&-\tau\leq\theta<0\\
      A_0(\eta)x(0)+A_\tau(\eta)x(-\tau)&,&\theta=0,
    \end{array}
\right.
\end{equation*}
where  $A_0(\eta)=\frac{\partial}{\partial x}f(x^\pm,x^\pm,\eta)$, $A_\tau(\eta)=\frac{\partial}{\partial x_\tau}f(x^\pm,x^\pm,\eta)$,  
$f$ is as in \eqref{eq:RFDE} and defined in \eqref{e.f}.
  
The characteristic equation for $A(\eta)$ is obtained by looking for nontrivial solution of the form $e^{\lambda t}c$ where $c\in \R^{2N}$ is a constant vector. Then $A(\eta)$  has an eigenvalue $\lambda$ with eigenfunction  $x(\theta) = e^{\lambda \theta}c$ if and only if
\begin{equation}
  \text{det}(\triangle(\lambda,\tau,\eta)):=\text{det}(\lambda \text{Id}-L(\eta,\tau))=0,
\end{equation}
where
\begin{equation}
  \label{eq:charac_matrix_def}
\triangle(\lambda,\tau,\eta):=\lambda\text{Id}-L(\eta,\tau)
\end{equation}
is the \textit{characteristic matrix} and  
\begin{equation}
  \label{eq:L(tau)_definition}
  L(\eta,\tau):=A_0(\eta)+A_\tau(\eta)e^{-\lambda\tau} \in \Mat(n)
\end{equation}
with $n=2N$,
see ~\cite{Hale1977, Wu1998}. We define the \textit{transcendental characteristic function} associated to 
$A(\eta)$ as
\begin{equation*}
  \label{eq:Transcendental_function}
  P(\lambda,\tau,\eta):=\text{det}(\triangle(\lambda,\tau,\eta)).
\end{equation*}
Since $F$ in  \eqref{e.F} is $\S_N$-equivariant, the matrix $\triangle(\lambda,\tau)$ is also $\S_N$-equivariant~\cite{Ruan2006}.
Thus $L$ from~\eqref{eq:L(tau)_definition} can be  decomposed as
\begin{equation*}
  L\cong\left(
    \begin{array}{cc}
     L_1 &  \\
          & L_{N-1} 
    \end{array}
\right),
\end{equation*}
and
\begin{equation*}
  L_{N-1}=\left(
          \begin{array}{cccc}
      L_2 &&&  \\
                  & L_2\\
                  &&\ddots\\
                  &&& L_2
    \end{array}
 \right),
\end{equation*}
$L_1,L_2\in\Mat(2,2)$. 
 Computing $L$ from~\eqref{e.f} we get:
 \marginnote{
\begin{align*}
& \frac{d }{\d \epsilon} \sin(x_{1,\tau}^{(j)}+ \epsilon v_{1,\tau}^{(j)}-
x_{1}^{(i)}- \epsilon v_{1}^{(i)}) \\
& = v_{1,\tau}^{(j)}- v_{1}^{(i)}=
c_1^{(j)} e^{-\lambda t} -c_1^{(i)} 
\end{align*}
 }
\begin{equation*}
  \label{eq:full_phase_linear_operator}
  Lx=\left(
    \begin{array}{c}
      x_2^{(1)}\\
      K \mu(-1+\cos(2\phi^\pm))x_1^{(1)}-\mu x_2^{(1)}+\dfrac{K\mu}{N-1}(1+\cos(2\phi^\pm))e^{-\lambda\tau}\displaystyle\sum_{\stackrel[j\neq 1]{j=1}{}}^Nx_1^{(j)}\\
\vdots\\
      x_2^{(N)}\\
     K \mu(-1+\cos(2\phi^\pm))x_1^{(N)}-\mu x_2^{(N)}+\dfrac{K\mu}{N-1}(1+\cos(2\phi^\pm))e^{-\lambda\tau}\displaystyle\sum_{\stackrel[j\neq N]{j=1}{}}^Nx_1^{(j)}
    \end{array}
\right) 
\end{equation*}
Hence, the characteristic matrix $\triangle(\lambda,\tau)=\lambda I_{2N}-L\in\Mat(2N)$ from ~\eqref{eq:charac_matrix_def}  has the  form
\begin{equation*}
  \triangle(\lambda,\tau)=\left(
    \begin{array}{cccc}
      m_{\lambda}&m_r&\cdots&m_r\\
      m_r&m_{\lambda}&\cdots&m_r\\
      \vdots&\vdots&\ddots&\vdots\\
      m_r&m_r&\cdots&m_{\lambda}
    \end{array}
\right)
\end{equation*}
where the blocks $m_{\lambda}$ and $m_r\in\text{Mat}(2)$ are
\begin{equation*}
  \label{eq:Nn_ml_mb_matrices}
  m_{\lambda}=\left(
    \begin{array}{cc}
      \lambda&-1\\
      q&\lambda+\mu
    \end{array}
\right),~~~m_r=\left(
  \begin{array}{cc}
    0&0\\r&0
  \end{array}
\right),
\end{equation*}
with  
\begin{equation}
\label{e.qr}
q=K\mu(1-\cos(2\phi^\pm)),\quad r=-\dfrac{K\mu}{N-1}(1+\cos(2\phi^\pm))e^{-\lambda\tau}. 
\end{equation}
Let $\Z_N$  be the cyclic group of order $N$ which is generated by the transformation $\zeta$ that sends $\phi_j$
to $\phi_{(j+1)\mod N}$, $j=0,\ldots, N-1$. Then each  (complex) irreducible representation of $\Z_N$, such that
$\zeta$ acts as $\e^{2\pi \i j/N}$, appears exactly once in the permutation representation of $\S_N$ on $\R^N$,
hence the $j$-th isotypic component $V_j$ of $\Z_N$ on $\R^{2N}$ is two-dimensional and spanned by the row vectors of the matrix
\begin{equation}\label{e.Wj}
 W_j=\frac{1}{\sqrt{N}}\left(\begin{array}{ccccccc}
    \lambda_{0j}&0&\lambda_{1j}&0&\ldots&\lambda_{(N-1)j}&0\\\\
    0&\lambda_{0j}&0&\lambda_{1j}&\ldots&0&\lambda_{(N-1)j}
  \end{array}\right),
\end{equation}
 where $\lambda_{kj}=\lambda_{(k\cdot j)\,\text{mod}\,N}=\e^{\i\pi((k\cdot j)\,\text{mod}\,N)/N}$
and the row vectors of $W_0$ span $\Fix(\S_N)$; then the restriction of the characteristic matrix   $\triangle(\lambda,\tau)$ to $V_j$ is 
\begin{align*}
\begin{array}{l}
  \triangle(\lambda,\tau)|_{V_j}= \overline{W_j}\triangle(\lambda,\tau)W_j^T\\\\
  = \dfrac{1}{N}
   (  \overline{\lambda}_{0j}I_2,\ldots,\overline{\lambda}_{(N-1)j}I_2  )\left(
    \begin{array}{cccc}
      m_{\lambda}&m_r&\cdots&m_r\\
      m_r&m_{\lambda}&\cdots&m_r\\
      \vdots&\vdots&\ddots&\vdots\\
      m_r&m_r&\cdots&m_{\lambda}
    \end{array}
    \right)\left(\begin{array}{c}\lambda_{0j}I_2\\\vdots\\\lambda_{(N-1)j}I_2\end{array}\right)\\\\
=\dfrac{1}{N}(\overline{\lambda}_{0j}m_{\lambda}+m_r\displaystyle\sum_{\stackrel[k\neq0]{k=0}{}}^{N-1}\overline{\lambda}_{kj},\ldots,\overline{\lambda}_{(N-1)j}m_{\lambda}+m_r\displaystyle\sum_{\stackrel[k\neq N-1]{k=0}{}}^{N-1}\overline{\lambda}_{kj} 
   )\left(\begin{array}{c}\lambda_{0j}I_2\\\vdots\\\lambda_{(N-1)j}I_2\end{array}\right)\\\\
= \dfrac{1}{N}\left(\displaystyle\sum_{\ell=0}^{N-1}\left(m_{\lambda}\lambda_{\ell j}+m_r\displaystyle\sum_{\stackrel[k\neq \ell]{k=0}{}}^{N-1}\lambda_{kj} \right)\overline{\lambda}_{\ell j}I_2 \right)\\\\
= \dfrac{1}{N}\left(\displaystyle\sum_{\ell=0}^{N-1}m_{\lambda}\lambda_{\ell j} \overline{\lambda}_{\ell j}+m_r\displaystyle\sum_{\ell=0}^{N-1}\displaystyle\sum_{\stackrel[k\neq \ell]{k=0}{}}^{N-1}\lambda_{kj}\overline{\lambda}_{\ell j} \right)\\\\
=m_{\lambda}+\dfrac{1}{N}m_r\displaystyle\sum_{\ell=0}^{N-1}\displaystyle\sum_{\stackrel[k\neq
    \ell]{k=0}{}}^{N-1}\lambda_{kj}\overline{\lambda}_{\ell j}.
\end{array}
\end{align*}
Moreover,
\begin{equation*}
\begin{array}{l}
\text{if}~j=0,~~~\displaystyle\sum_{\ell=0}^{N-1}\displaystyle\sum_{\stackrel[k\neq \ell]{k=0}{}}^{N-1}\lambda_{kj}\overline{\lambda}_{\ell j}=N(N-1)\\\\
\text{if}~j\neq 0,~~~\displaystyle\sum_{\ell=0}^{N-1}\displaystyle\sum_{\stackrel[k\neq \ell]{k=0}{}}^{N-1}\lambda_{kj}\overline{\lambda}_{\ell j}=
-\displaystyle\sum_{\ell=0}^{N-1}\lambda_{\ell j}\overline{\lambda}_{\ell j} =-N.
\end{array}
\end{equation*}
Therefore
\begin{equation*}
  \label{eq:Nn_Delta_red}
  \triangle(\lambda,\tau)|_{V_j}=\left\{
    \begin{array}{ll}
      m_{\lambda}+(N-1)m_r,&j=0,\\\\
      m_{\lambda}-m_r,&j=1,\ldots,N-1.
    \end{array}
\right.
\end{equation*}
The characteristic matrix decomposition is
\begin{equation*}
  \label{eq:Nn_Delta_4}
\triangle(\lambda,\tau)=\text{diag}\left(
      \triangle(\lambda,\tau)|_{\Fix(\S_N)},\triangle(\lambda,\tau)|_{V_1},\ldots,\triangle(\lambda,\tau)|_{V_{N-1}}
\right),
\end{equation*}
where $ {\Fix(\S_N)}= {V_0}$.

The characteristic function $P(\lambda,\tau)$ defined in~\eqref{eq:Transcendental_function} becomes
\begin{equation*}
P(\lambda,\tau)=\det(\triangle(\lambda,\tau)|_{\Fix(\S_N)})\prod_{j=1}^{N-1}\det(\triangle(\lambda,\tau)|_{V_j}),
\end{equation*}
or
\begin{equation}
  \label{eq:Nn_charact_functions_2}
  P(\lambda,\tau)=\text{det}(m_{\lambda}+(N-1)m_r)(\text{det}(m_{\lambda}-m_r ))^{N-1},
\end{equation}
and using~\eqref{eq:Nn_ml_mb_matrices} we obtain
\begin{equation}
  \label{eq:full_phase_P}
  \begin{array}{rcl}
    P_{\Fix(\S_N)}(\lambda,\tau)&=&\det(\triangle(\lambda,\tau)|_{\Fix(\S_N)}) \\
    &=&\lambda^2+\mu\lambda+q+(N-1)r\\\\
    P_U(\lambda,\tau)&=&\det(\triangle(\lambda,\tau)|_{V_j}) ,~~j\neq 0,\\
    &=&\lambda^2+\mu\lambda+q -r ,
  \end{array}
\end{equation} 
 where $q$ and $r$ are as in \eqref{e.qr} and $U$ as in \eqref{e.U}.

For the equilibria in~\eqref{eq:full_phase_equilibria_points} we have
\begin{equation*}
  \cos(2\phi^\pm)=\pm \dfrac{1}{K}\sqrt{K^2-\omega_M^2},
\end{equation*}
with $K\geq\omega_M$ in order to keep $\phi^\pm\in\mathbb{R}$. By scaling $\widetilde{K}=K/\omega_M$, $\tilde{\mu}=\mu/\omega_M$, $\tilde{\lambda}=\lambda/\omega_M$, and $\tilde{\tau}=\omega_M\tau$, and removing the tilde in the variables we obtain the normalized equilibria:
\begin{equation}
  \label{eq:full_phase_phi**}
  \cos(2\phi^\pm)= 
      \pm\sqrt{1-\dfrac{1}{K^2}},
\end{equation}
with
\begin{equation}
\label{eq:full_phase_K_condition}
 K\geq1.
\end{equation}


\subsection{Symmetry-preserving  bifurcations}
\label{subsec:Bif_Pfix_full_phase_model}
In the next two sections bifurcations in the two isotypic components
$\Fix(\S_N)$ and $U\oplus U$ of $\S_N$ found previously
are analyzed, conditions for the existence of eigenvalues $\lambda=\pm \i\omega$ with $\omega\in\R^+$ are  given in terms of the parameters $K,\mu,\tau\in\R^+$, and $N\in\N>1$  and the critical time delays $\tau$ leading to bifurcation are computed. When $\tau=0$ the transcendental characteristic functions in~\eqref{eq:full_phase_P} become ordinary characteristic polynomials with two roots each. Since we are interested in analyzing the influence of the time-delay between the nodes in the network it is important to know whether the system is stable or not at $\tau=0$. If it is, we would like to determine, if there exists some $\tau\in\R^+$ such that a finite number of roots cross the imaginary axis at $\lambda=i \omega$ from the left to the right switching stability with $d\lambda/d\tau|_{\lambda=\i\omega}\neq0$. For this analysis we use the $S_n$ map which 
we discuss in section~\ref{subsec:Sn_map}. If some roots are unstable at the equilibrium at $\tau=0$ we look for some $\tau\in\R^+$ such that all unstable roots (always a finite number) cross from the right to the left  at $\tau$ switching stability from unstable to stable; this task is addressed using the Lambert W function see~\cite{Asl2003, Corless1996, Mathews2004, Wang2008}.
We start by analyzing bifurcations in $\Fix(\S_N)$. Note that $\S_N$-symmetry implies that \eqref{e.F} maps $\Fix_\kX(\S_N)$ to itself, which means that we can restrict \eqref{e.fullPhase} to the  subspace $\Fix(\S_N)$~\cite{Golubitsky1988}.

Note that bifurcations in $\Fix(\S_N)$, which we study first, preserve the $\S_N$-symmetry, so bifurcating periodic orbits are synchronized, i.e., satisfy, $\phi_i=\phi_j$ for all $i,j$. In section \ref{subsec:Bif_Pj_full_phase_model} we will study bifurcations in the other blocks which are symmetry-breaking, i.e., bifurcating periodic orbits are not fully-synchronized, for more details
see Section \ref{subsec:spatio-temp_symm}.


\subsubsection{Roots of the characteristic function $P_{\Fix(\S_N)}(\lambda,\tau)$ at $\tau=0$ and as $\tau\to\infty$}
\label{subsubsec:roots_Pfix_tau0_tau_infty}
In the fixed-point space in equation~\eqref{eq:full_phase_P} when $\tau=0$ we have two roots
\begin{equation}
  \lambda_{\pm}=-\frac{1}{2}\mu\pm\frac{1}{2}\left(\mu^2+8K\mu\cos(2\phi^\pm) \right)^{1/2};
\end{equation}
here we have two cases corresponding to $\phi^\pm$
from~\eqref{eq:full_phase_phi**} 
\begin{equation}
  \label{eq:full_phase_charc_fun_tau0_a}  
  \lambda_{\pm}=-\frac{1}{2}\mu\pm\frac{1}{2}\left(\mu^2\pm8\mu\sqrt{K^2-1} \right)^{1/2}, 
\end{equation}
and remembering $K\geq1$, see~\eqref{eq:full_phase_K_condition}, we obtain that
\begin{itemize}
\item If $K>1$, there is an unstable root for $\phi^+$, and both roots are stable for $\phi^-$.
\item If $K=1$, there is a constant root at $\lambda=0$ and another one at $\lambda=-\mu$, for the unique equilibrium $\phi^+=\phi^-$.
\end{itemize}
In the limit when $\tau\to\infty$ in equation~\eqref{eq:full_phase_P} for both equilibria $\phi^\pm$, assuming that $\Re(\lambda)>0$, we obtain
\begin{equation}
  \label{eq:full_phase_roots_tau_infty}
\lambda_{\pm}=-\frac{1}{2}\mu\pm\frac{1}{2}\left(\mu^2-4K\mu\left(1\mp\sqrt{1-\dfrac{1}{K^2}}\right) \right)^{1/2},
\end{equation}
\marginnote{\tiny
Since $\sqrt{1- \frac{1}{K^2}} \in [0,1)$ we see that
$\left(1\mp\sqrt{1-\dfrac{1}{K^2}}\right)   \in (0,2]$, so the whole
square root is either real and in $[0,\mu)$ or purely imaginary.
}
but these roots are not in the right-side of the complex plane neither for $\phi^+$ nor $\phi^-$, which is a contradiction, therefore at $\tau\to\infty$ both equilibria are spectrally stable in $\Fix(\S_N)$.


\subsubsection{The $S_n$ map}
\label{subsec:Sn_map}
In~\cite{Beretta2002}    a criterium is presented  to find imaginary roots for a transcendental function of the form
\begin{equation}
  \label{eq:Sn_map_P}
  P(\lambda,\tau)=R(\lambda,\tau)+S(\lambda,\tau)e^{-\lambda\tau},
\end{equation}
where
\begin{equation}
  \label{eq:Sn_map_R_S}
  R(\lambda,\tau)=\sum_{k=0}^nr_k(\tau)\lambda^k,~~~S(\lambda,\tau)=\sum_{k=0}^ms_k(\tau)\lambda^k.
\end{equation}
In~\eqref{eq:Sn_map_R_S}, $n,m\in\N_0$, $n>m$, and $r_k(\cdot),~s_k(\cdot):\R^+_0\rightarrow\R$ are continuous and differentiable functions of $\tau$. We shall describe  the method briefly and then apply it to the full-phase model, and other
models subsequently.

We are looking for roots $\lambda=\pm \i\omega$ of $P(\lambda,\tau)$ from~\eqref{eq:Sn_map_P}, with $\omega\in\R^+_0$. Since the roots appear in complex conjugate pairs, we only need to  look   for roots with $\omega\geq 0$. Substituting $\lambda=\i\omega$ into~\eqref{eq:Sn_map_P} we have
\marginnote{
\tiny
\begin{align*} e^{-i\omega \tau} &= -\frac{R}{S} = -\frac{R \bar S}{|S|^2}\\
&\frac{-(R_R+i R_I)(S_R-iS_I)}{|S|^2}
\end{align*}
}
\begin{equation}
  \label{eq:Sn_map_sin_cos}
  \begin{array}{lll}
    \sin(\omega\tau)&=&\dfrac{R_IS_R- S_IR_R}{|S|^2}\\\\
    \cos(\omega\tau)&=&-\dfrac{S_IR_I+S_RR_R}{|S|^2}
  \end{array},~~~|S|\neq0,
\end{equation}
where $S_R$, $S_I$, $R_R$, and $R_I$ stand for the real part and the imaginary part of $S(\i\omega,\tau)$ and $R(\i\omega,\tau)$ respectively.

On the other hand, we can eliminate the exponential term  in \eqref{eq:Sn_map_P} and define the polynomial  equation
in $\omega$
\marginnote{
\tiny
\begin{align}
  F(\omega,\tau)&:=R(i\omega,\tau)R(-i\omega,\tau)-\nonumber\\
                &S(i\omega,\tau)S(-i\omega,\tau)=0\nonumber
\end{align}
\normalsize
}
\begin{equation}
\label{eq:Sn_map_F}
  F(\omega,\tau):=|R(\i\omega,\tau)|^2-|S(\i\omega,\tau)|^2=0. 
\end{equation}
\begin{definition}{}
Let $x,y\in\R$, $r\in\R^+$, and let $\theta\in(-\pi,\pi]$ satisfy $x=r\cos\theta$ and $y=r\sin\theta$. We define the argument of $(x,y)$ as $\arg\left(\cdot\right):\mathbb{R}^2\backslash\{(0,0)\}\to(-\pi,\pi]$, such that $\arg\left(x,y\right)=\theta$. This function is the extension of the trigonometrical function $\arctan(y/x)$ where $\arctan:\R\to(-\pi/2,\pi/2)$.
\end{definition}
Now, given $\tau\in\R^+$ we can compute possible values of
$\omega=\omega(\tau)$ as roots of the polynomial $F$
from~\eqref{eq:Sn_map_F}. Since $\sin(\omega\tau)$ and
$\cos(\omega\tau)$ in ~\eqref{eq:Sn_map_sin_cos} are both functions of
$\omega(\tau)$ and $\tau$, we can calculate the argument
$\theta(\tau)=\omega\tau+2n\pi$, for $n\in\Z$ using~\eqref{eq:Sn_map_sin_cos}  as
\begin{equation}
\theta(\tau)=\arg\left(-S_IR_I-S_RR_R, R_IS_R-S_IR_R \right).
\end{equation}
Then we define the map $\tau_n:\R^+_0\to\R$ as
\begin{equation}
  \label{eq:2n_tau_map}
  \tau_n(\tau):=\frac{\theta(\tau)+2n\pi}{\omega(\tau)}.
\end{equation}
If $\tau_n(\tau)=\tau$, then $\tau=\tau^*$ is a bifurcation time-delay, and $\i\omega(\tau^*)$ is an  imaginary root of \eqref{eq:Sn_map_P}; this can be formally expressed by the map
\begin{equation*}
  \label{eq:2n_Sn_map}
  S_n:=\tau-\tau_n(\tau),
\end{equation*}
whose zeros are the critical bifurcation time delays for equation~\eqref{eq:Sn_map_P}.

Now, we need to know in which direction the roots found above cross
the imaginary axis when $\tau$ is varied, if they go from stable to unstable or from unstable to stable in the complex plane. We need to calculate
\begin{equation}
  \label{eq:2n_delta}
\delta(\omega(\tau^*)):=\Re\left(\frac{d\lambda}{d\tau}\bigg|_{\lambda=\i\omega(\tau^*)}\right)=\Re\left(-\frac{dP}{d\tau}\left(\frac{dP}{d\lambda}\right)^{-1}_{\lambda=\i\omega(\tau^*)} \right);
\end{equation}
from the definition of $P(\lambda,\tau)$ in ~\eqref{eq:Sn_map_P} we have
\begin{equation}
  \frac{d\lambda}{d\tau}\bigg|_{\lambda=\i\omega(\tau^*)}=\frac{e^{-\i\omega\tau}(\i\omega S(\i\omega,\tau)-S'_{\tau}(\i\omega,\tau))-R'_{\tau}(\i\omega,\tau)}{R'_{\lambda}(\i\omega,\tau)+e^{-\i\omega\tau}(S'_{\lambda}(\i\omega,\tau)-\tau S(\i\omega,\tau))},
\end{equation}
where $R'_{\lambda}$ means derivative of $R(\lambda,\tau)$ with respect to $\lambda$, $R'_{\tau}$ derivative of $R(\lambda,\tau)$ with respect to $\tau$ and the same for $S(\lambda,\tau)$. Then we have
\begin{equation}
  \label{eq:Sn_delta_abcd}
  \delta(\omega(\tau^*))=\frac{AC+BD}{C^2+D^2},
\end{equation}
where
\begin{equation}
  \label{eq:Sn_map_abcd}
 \begin{array}{l}
  A=\Re(e^{-\i\omega\tau}(\i\omega S(\i\omega,\tau)-S'_{\tau}(\i\omega,\tau))-R'_{\tau}(\i\omega,\tau))
\\B=\Im(e^{-\i\omega\tau}(\i\omega S(\i\omega,\tau)-S'_{\tau}(\i\omega,\tau))-R'_{\tau}(\i\omega,\tau))
\\C=\Re(R'_{\lambda}(\i\omega,\tau)+e^{-\i\omega\tau}(S'_{\lambda}(\i\omega,\tau)-\tau S(\i\omega,\tau)))
\\D=\Im(R'_{\lambda}(\i\omega,\tau)+e^{-\i\omega\tau}(S'_{\lambda}(\i\omega,\tau)-\tau S(\i\omega,\tau))).
 \end{array}
\end{equation}
If $\delta(\omega(\tau^*))>0$ the root crosses from the left to the right (stable to unstable), and if $\delta(\omega(\tau^*))<0$ the root crosses in the opposite direction. It is important to note that condition $\delta(\omega(\tau^*))\neq0$, called transversality condition, is necessary for Hopf bifurcation to occur~\cite{Hassard1981}.%


\subsubsection{Conditions for the existence of  symmetry-preserving bifurcations}
\label{subsub:Condition_existence_bif_in_Pfix}
Since our aim is to analyze bifurcations in $\Fix(\S_N)$, we check the necessary conditions for the existence of roots $\lambda=\pm \i\omega$, $\omega\in\R^+$ given by the polynomial $F(\omega)$ from  \eqref{eq:Sn_map_F}.
From~\eqref{eq:Sn_map_P},~\eqref{eq:Sn_map_R_S} with $n=2$ and $m=0$  and \eqref{eq:full_phase_P} we have
\begin{equation}
  \label{eq:full_phase_R_S_pol}
\begin{array}{l}
  R(\lambda)=\lambda^2+\mu\lambda+K\mu(1-\cos(2\phi^\pm))
\\S(\lambda)=-K\mu(1+\cos(2\phi^\pm)),
\end{array}
\end{equation}
then $F(\omega)$ become
\marginnote{
\tiny
\begin{align*}
  R(\i \omega) &= -\omega^2 + \mu \i \omega + K \mu (-)  \Longrightarrow\\
 R(\i \omega)R(-\i \omega) &=\omega^4 + K^2 \mu^2 (-)^2 -2 \omega^2K \mu (-) \\
&  +  \mu^2  \omega^2\\
S(\i \omega)S(-\i \omega) &=K^2 \mu^2(+)^2\Longrightarrow\\
F &= \omega^4 + K^2 \mu^2( (-)^2 - (+)^2 ) \\
&+  \omega^2( \mu^2 -2 K \mu (-) )
\end{align*}
\normalsize
}
\begin{equation}
  \label{eq:full_phase_exist_cond}
  F(\omega)=\omega^4+\left(\mu^2-2K\mu\left(1-\cos(2\phi^*) \right)  \right)\omega^2-4K^2\mu^2\cos(2\phi^*),
\end{equation}
and
\begin{equation}
  \label{eq:full_phase_w}
\begin{array}{rcl}
 \omega_{\pm}^2&=&-\dfrac{1}{2}\left(\mu^2-2K\mu\left(1-\cos(2\phi^*)\right) \right)\\
&&\pm\dfrac{1}{2}\left[\left(\mu^2-2K\mu\left(1-\cos(2\phi^*)\right) \right)^2+16K^2\mu^2\cos(2\phi^*)\right]^{1/2},
\end{array}
\end{equation}
where $\phi^* = \phi^\pm$.

For the sake of notation we write
\begin{equation}
  \label{eq:w_pm_short_form}
  \omega_{\pm}^2=-\frac{1}{2}b\pm\frac{1}{2}\sqrt{b^2-4c},
\end{equation}
where
\begin{equation}\label{e.bc}
  \begin{array}{rcl}
    b&=&\mu^2-2K\mu\left(1-\cos(2\phi^*)\right) \\
    c&=&-4K^2\mu^2\cos(2\phi^*),
  \end{array}
\end{equation}
with $\cos(2\phi^*)=\pm\sqrt{1-(1/K^2)}$.
\begin{lemma}
\label{l.conditions_bif_Xj}
 A necessary condition for the existence of $\omega_{\pm}\in\R$ is
\begin{equation}
\label{eq:full_phase_Pj_general_exist_cond_a}
b^2-4c\geq0.
\end{equation}
Moreover,  if \eqref{eq:full_phase_Pj_general_exist_cond_a} holds then:
\begin{enumerate}
\item  If $b\geq0$
\begin{enumerate}
 \item  If $c\leq0$ then $\omega_+\in\R^+_0$ and, if $bc \neq 0$, then $\omega_-\in\C\setminus\R$.
 \item If $c>0$ then $\omega_{\pm}\in\C\setminus\R$.
\end{enumerate}
\item  If $b<0$ 
  \begin{enumerate}
  \item  If $c\leq0$ then $\omega_{+}\in\R^+$, $\omega_{-}\in\C\setminus\R$.
  \item  If $c>0$ then $\omega_{\pm}\in\R^+$.
\end{enumerate}
\end{enumerate}
\end{lemma}
Provided $\omega_{+} \in \R$ or $\omega_- \in \R$, we can find the critical time-delay $\tau\in\R^+$ such that
$\i \omega$ is a root of $P_{\Fix(\S_N)}$ using the $S_n$ map in section~\ref{subsec:Sn_map}, thus, from~\eqref{eq:full_phase_R_S_pol} and~\eqref{eq:2n_tau_map} we have
\marginnote{
\reversemarginpar
\tiny
\begin{align}
 \sin(\omega_{\pm}\tau)&=\dfrac{R_IS_R-S_IR_R}{|S|^2}\nonumber\\
                      &=\dfrac{R_IS_R}{|S_R|^2}\nonumber\\
                      &=-\dfrac{\omega_{\pm}}{K(1+\cos(2x_1^*))}\nonumber\\
 \cos(\omega_{\pm}\tau)&=-\dfrac{S_IR_I+S_RR_R}{|S|^2}\nonumber\\
                      &=-\dfrac{R_RS_R}{|S_R|^2}\nonumber\\
                      &=\dfrac{-\omega_{\pm}^2+K\mu(1-\cos(2x_1^*))}{K\mu(1+\cos(2x_1^*))}.\nonumber
\end{align}
\normalsize
}
\begin{equation}
  \label{eq:full_phase_sin_cos_Pfix}
  \begin{array}{rcl}
    \sin(\omega_{\pm}\tau)&=&-\dfrac{\omega_{\pm}}{K(1+\cos(2\phi^*))}\\\\
    \cos(\omega_{\pm}\tau)&=&\dfrac{-\omega_{\pm}^2+K\mu(1-\cos(2\phi^*))}{K\mu(1+\cos(2\phi^*))},
  \end{array}
\end{equation}
and 
\begin{equation}
  \label{eq:full_phase_tau_general}
  \tau_{\pm}(\omega_{\pm};K,\mu,n)=\frac{1}{\omega_{\pm}}\left[\arg\left( \cos(\omega_{\pm}\tau), \sin(\omega_{\pm}\tau)  \right)+ 2n\pi \right],~~~n\in\mathbb{Z}.
\end{equation}
Here we want to stress  that $\tau_{\pm}$ does not depend on $\tau$, see \eqref{eq:full_phase_sin_cos_Pfix}. In what follows we will write $\tau_{\pm}(n) = \tau(\omega_{\pm},n)$  to emphasize the dependence on $\omega_{\pm}$ or $n$ respectively according to need.

The direction in which the roots cross the imaginary axis, if there are any, can be obtained by looking at the sign of $\delta(\omega)$ defined in \eqref{eq:Sn_delta_abcd},
where, due to \eqref{eq:full_phase_R_S_pol} and \eqref{eq:full_phase_sin_cos_Pfix},
 the constants from \eqref{eq:Sn_map_abcd} are
\marginnote{
\tiny
\begin{align}
   A=&~\Re(e^{-i\omega\tau}(i\omega S(i\omega,\tau)\nonumber
\\   &~-S'_{\tau}(i\omega,\tau))-R'_{\tau}(i\omega,\tau))\nonumber
\\  =&~\Re\left[\left(\cos(\omega_{\pm}\tau)-i\sin(\omega_{\pm}\tau)\right)(i\omega_{\pm})\right.\nonumber
\\   &\left.\times(-K\mu(1+\cos(2x_1^*))) \right]\nonumber
\\  =&~-K\mu\omega_{\pm}(1+\cos(2x_1^*))\sin(\omega_{\pm}\tau)\nonumber
\\  =&~-K\mu\omega_{\pm}(1+\cos(2x_1^*))\nonumber
\\   &~\times\left(\dfrac{-\omega_{\pm}}{K(1+\cos(2x_1^*))}\right)\nonumber
\\  =&~\mu\omega_{\pm}^2\nonumber
\\ B=&~\Im(e^{-i\omega\tau}(i\omega S(i\omega,\tau)\nonumber
\\   &~-S'_{\tau}(i\omega,\tau))-R'_{\tau}(i\omega,\tau))\nonumber
\\  =&~-K\mu\omega_{\pm}(1+\cos(2x_1^*))\cos(\omega_{\pm}\tau)\nonumber
\\  =&~-K\mu\omega_{\pm}(1+\cos(2x_1^*))\nonumber
\\   &~\times\left(\dfrac{-\omega_{\pm}^2+K\mu(1-\cos(2x_1^*))}{K\mu(1+\cos(2x_1^*))}\right)\nonumber
\\  =&~\omega_{\pm}^3-K\mu\omega_{\pm}(1-\cos(2x_1^*))\nonumber
\\ C=&~\Re(R'_{\lambda}(i\omega,\tau)\nonumber
\\   &~+e^{-i\omega\tau}(S'_{\lambda}(i\omega,\tau)-\tau S(i\omega,\tau)))\nonumber
\\  =&~\Re\left[(2i\omega_{\pm}+\mu)+(\cos(\omega_{\pm}\tau)\right.\nonumber
\\   &\left.-i\sin(\omega_{\pm}\tau))(\tau K\mu(1+\cos(2x_1^*))) \right]\nonumber
\\  =&~\mu(1+K\tau_{\pm}(1+\cos(2x_1^*))\cos(\omega_{\pm}\tau))\nonumber
\\  =&~\mu+\mu K\tau_{\pm}(1+\cos(2x_1^*))\nonumber
\\   &~\times\left(\dfrac{-\omega_{\pm}^2+K\mu(1-\cos(2x_1^*))}{K\mu(1+\cos(2x_1^*))}\right)\nonumber
\\  =&~\mu-\tau_{\pm}\omega_{\pm}^2+\tau_{\pm}K\mu(1-\cos(2x_1^*))\nonumber
\\ D=&~\Im(R'_{\lambda}(i\omega,\tau)\nonumber
\\   &~+e^{-i\omega\tau}(S'_{\lambda}(i\omega,\tau)-\tau S(i\omega,\tau)))\nonumber
\\  =&~2\omega_{\pm}-\mu K\tau(1+\cos(2x_1^*))\sin(\omega_{\pm}\tau)\nonumber
\\  =&~2\omega_{\pm}-\mu K\tau(1+\cos(2x_1^*))\nonumber
\\   &~\times\left(\dfrac{-\omega_{\pm}}{K(1+\cos(2x_1^*))}\right)\nonumber
\\  =&~\omega_{\pm}(2+\mu\tau),\nonumber
\\&\nonumber
\\&\nonumber
\\AC+&BD=\nonumber
\\     =&~\mu\omega_{\pm}^2\left(-\tau_{\pm}\omega_{\pm}^2\right.\nonumber
\\        &~\left.+\mu(1+\tau_{\pm}K(1-\cos(2x_1^*)))\right)\nonumber
\\      &~+\left(\omega_{\pm}^3-K\mu\omega_{\pm}(1-\cos(2x_1^*)) \right)\nonumber
\\      &~(2\omega_{\pm}+\mu\tau\omega_{\pm})\nonumber
\\     =&~-\mu\tau\omega_{\pm}^4\nonumber
\\      &~+\mu^2(1+\tau K(1-\cos(2x_1^*)))\omega_{\pm}^2\nonumber
\\        &~+\omega_{\pm}^4(2+\mu\tau)\nonumber
\\      &~-K\mu(1-\cos(2x_1^*))(2+\mu\tau)\omega_{\pm}^2\nonumber
\\ =&\omega_{\pm}^2\left(-\mu \tau^2\omega^2 +\mu^2(1+\tau K(-)) \right.\nonumber
\\  &\left.+ 2\omega^2
+ \omega^2\mu\tau -K\mu(-)(2+\mu\tau) \right) \nonumber
\\     =&~\omega_{\pm}^2\left(2\omega_{\pm}^2+\right.\nonumber
\\      &~\left.\mu^2-2K\mu(1-\cos(2x_1^*))\right)\nonumber
\\     =~&\omega_{\pm}^2\left(2\omega_{\pm}^2+b\right)\nonumber 
\end{align}
\normalsize
}
\begin{align}
  \label{eq:full_phase_delta_abcd_general}
\begin{array}{l}
    A= \mu\omega_{\pm}^2\\
    B= \omega_{\pm}^3-K\mu\omega_{\pm}(1-\cos(2\phi^*))\\
    C= \mu-\tau_{\pm}\omega_{\pm}^2+\tau_{\pm}K\mu(1-\cos(2\phi^*))\\
    D= \omega_{\pm}(2+\mu\tau_{\pm}).
\end{array}
\end{align}
It is clear from \eqref{eq:Sn_delta_abcd} that the sign of $\delta$ depends on the numerator $AC+BD$, then using~\eqref{eq:full_phase_delta_abcd_general} we compute
\begin{equation}
  \label{eq:full_phase_AC+BD}
  \begin{array}{rcl}
  AC+BD&=&\omega_{\pm}^2\left(2\omega_{\pm}^2+b\right) 
  \end{array}
\end{equation}
but from equation~\eqref{eq:w_pm_short_form} we know that $b=-2\omega_{\pm}^2\pm\sqrt{b^2-4c}$, then substituting into~\eqref{eq:full_phase_AC+BD} we have
\begin{equation}
\label{eq:full_phase_delta_sign_a}
  AC+BD=\pm\omega_{\pm}^2\sqrt{b^2-4c},
\end{equation}
thus the sign of $\delta (\omega_{\pm}) $ is
\begin{equation}
  \label{eq:full_phase_AC+BD_b}
  \text{sgn}(\delta (\omega_{\pm}))=\left\{\begin{array}{l}
                                          +1~~~\text{for}~~~\omega_+\\
                                          -1~~~\text{for}~~~\omega_-\\       
                                          \end{array} \right..
\end{equation}


\subsubsection{Curves of symmetry-preserving bifurcations}
\label{subsub:Bifurcation_curves_Fix}
In this section we shall analyze the bifurcation curves in $\Fix(\S_N)$, from which fully synchronized
periodic orbits emanate, in three cases:
\begin{itemize}
\item When $K=1$ the following is valid for the unique equilibrium $\phi^*=\phi^+=\phi^-\in (-\pi,\pi]$. In this case the roots of the characteristic function $P_{\Fix(\S_N)}$ from \eqref{eq:full_phase_P} when $\tau=0$ are, by \eqref{eq:full_phase_charc_fun_tau0_a},
\begin{equation*}
  \lambda_+=0,~~~\lambda_-=-\mu.
\end{equation*}
For $\tau\neq 0$ the equation $F(\omega)=0$ in~\eqref{eq:full_phase_exist_cond}, which represents a necessary condition for the existence of roots at $\lambda=\pm \i\omega$, becomes, due to \eqref{eq:full_phase_phi**},
\begin{equation*}
  F(\omega)=\omega^4+\left(\mu^2-2\mu \right)\omega^2,
\end{equation*}
and here, except from the zero root which exists
for all $\tau$ due to a saddle node bifurcation at $K=1$, we have the following root
\begin{equation}
  \label{eq:full_phase_othe_roots}
  \omega =\pm\sqrt{2\mu-\mu^2}
\end{equation}
which is real if
\begin{equation}
  \label{eq:full_phase_other_cond}
  0<\mu<2.
\end{equation}
If ~\eqref{eq:full_phase_other_cond} does not hold, the roots of $P_{\Fix(\S_N)}$   remain in the left hand side of the complex plane with a constant root at zero, for all $\tau\in\R^+$. 

From~\eqref{eq:full_phase_sin_cos_Pfix}   we obtain
\marginnote{
\tiny
\begin{align}
  \sin(\omega_{+}\tau)&=~\dfrac{R_IS_R-S_IR_R}{|S|^2}\nonumber
\\                  &=~-\omega_{+}\nonumber
\\\cos(\omega_{+}\tau)&=~-\dfrac{S_IR_I+S_RR_R}{|S|^2}\nonumber
\\                    &=\dfrac{\mu-\omega_{+}^2}{\mu}\nonumber
\end{align}
\normalsize
}
\begin{equation}
  \label{eq:full_phase_K=wm_sin_cos}
  \begin{array}{rcl}
    \sin(\omega \tau)&=&-\omega \\
    \cos(\omega \tau)&=&\dfrac{\mu-\omega^2}{\mu}.
  \end{array}
\end{equation}
From~\eqref{eq:full_phase_K=wm_sin_cos} we compute $\tau$ as a function
of $\omega$ and $\mu$,
\begin{equation}
  \label{eq:full_phase_tau}
  \tau(\omega;\mu,n)=\frac{1}{\omega}\left[ 
  \arg\left(\mu-\omega^2,-\omega\mu\right)+2n\pi \right],~~~n\in\Z.
\end{equation}
We already know from \eqref{eq:full_phase_AC+BD_b} that the roots $\lambda=\pm\i\omega$ cross the imaginary axis from the left to the right. 
In figure~\ref{fig:mu1_x_tau_full_phase} the curves for $\tau(\omega;\mu,n)$ from \eqref{eq:full_phase_tau}  are shown 
for different values of $n$.

  The curve $\tau(n)=\tau(\omega;\mu,n)$ for $n=0$ determines the first root crossing from the left to the right. For each curve $\tau(n)$ with $n\in \N_0$ a new root crosses from the left-hand side to the right-hand side of the imaginary axis.
\begin{figure}[!htb]
  \centering
  \includegraphics[scale=0.4]{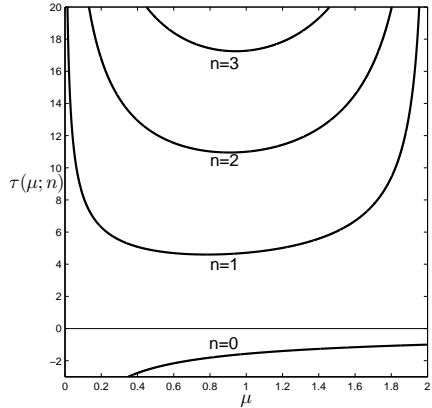}
  \caption{Symmetry-preserving bifurcation curves for the equilibrium $\phi^+=\phi^-$, with $K=1$ and $0<\mu<2$.}
  \label{fig:mu1_x_tau_full_phase}
\end{figure}
\begin{remark}
Although 
 the equilibrium  $\phi^+=  \phi^-$ is  spectrally stable in $\Fix(\S_N)$  (no
roots of $P_{\Fix(\S_N)}$ with positive real part) with $K=1$ at
$\tau=0$, see \eqref{eq:full_phase_charc_fun_tau0_a}, and the first bifurcation root 
appears on the curve $\tau(\mu;0)$, we cannot conclude  stability of the equilibrium in $\Fix(\S_N)$ below this curve due a constant zero root $\lambda=0$ caused by a saddle node bifurcation of the curve of $\S_N$-invariant equilibria
$\phi(K)$ given by $\phi^\pm(K)$ at $K=1$. 
Note that $K$  
represents the coupling strength between nodes and stability when this parameter
varies has been already studied in literature: in~\cite{Earl2003}, a
stability criterion for a
general coupling function is derived, and in~\cite{Yeung1999} the
stability of the Kuramoto model is studied;
an extensive review of these and other related results can be found
 in~\cite{Klinshov2013}.
\end{remark}

\item Next we analyze the unstable equilibrium $\phi^+$ when
 $K>1$. We are interested in any values of parameters $\mu$ and $\tau$ such that the roots in $P_{\Fix(\S_N)}(\lambda,\tau)$ become stable, i.e.,   in any $\mu,\tau\in\R^+$ for which we have $\max(\Re\lambda)<0$, for  all roots $\lambda$ of  $P_{\Fix(\S_N)}$. We shall use the Lambert W function, see~\cite{Asl2003, Corless1996, Mathews2004, Wang2008}, to find the rightmost root when $\mu$ and $\tau$ vary. 
The initial $\lambda_0$ guess needed in both Newton's and Halley's schemes used to calculate the rightmost root is found using the rightmost root in an auxiliary polynomial as proposed in~\cite{Wang2008, Wang2008a}, and for the following iterations the root found in the previous iteration is used as initial guess. Results of the numerical simulation with $\mu=\{0.1,~0.2,~0.4,~0.6,~0.8\}$ and $K=2$ are shown in figure~\ref{fig:fix_point_full-phase_rightmost_root}. As expected at $\tau=0$ the real part of the rightmost root is positive and increases monotonically with $\mu$, see~\eqref{eq:full_phase_charc_fun_tau0_a}. On the other hand, when $\tau$ grows the real part of the rightmost root tends to a non positive value as shown in~\eqref{eq:full_phase_roots_tau_infty} and predicted
in Section \ref{subsubsec:roots_Pfix_tau0_tau_infty}.
\begin{figure}[!htb]
  \centering
  \begin{subfigure}[t]{2.2in}
  \centering
  \includegraphics[scale=0.4]{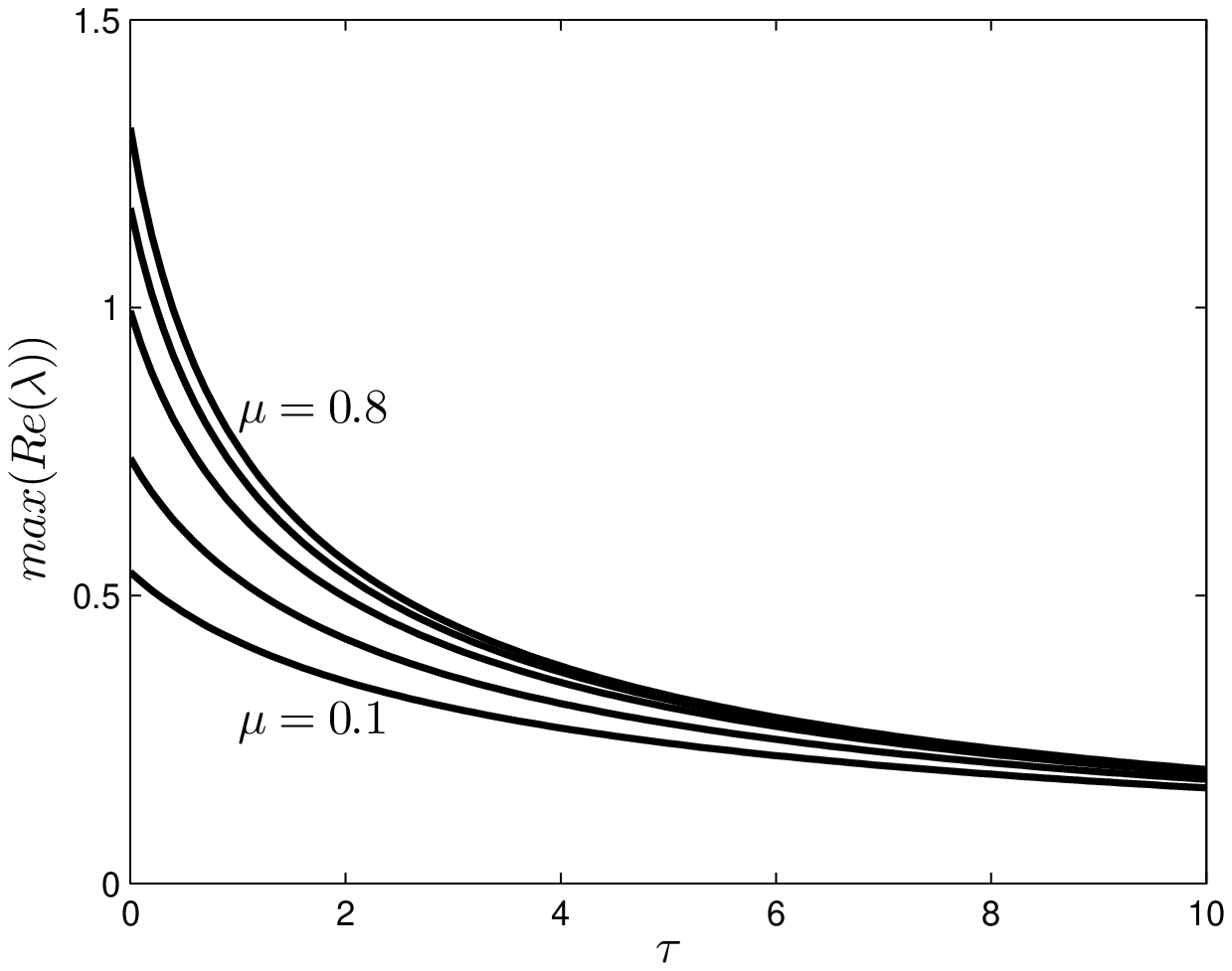}
  \end{subfigure}
  \quad
  \begin{subfigure}[t]{2.2in}  
  \centering
  \includegraphics[scale=0.4]{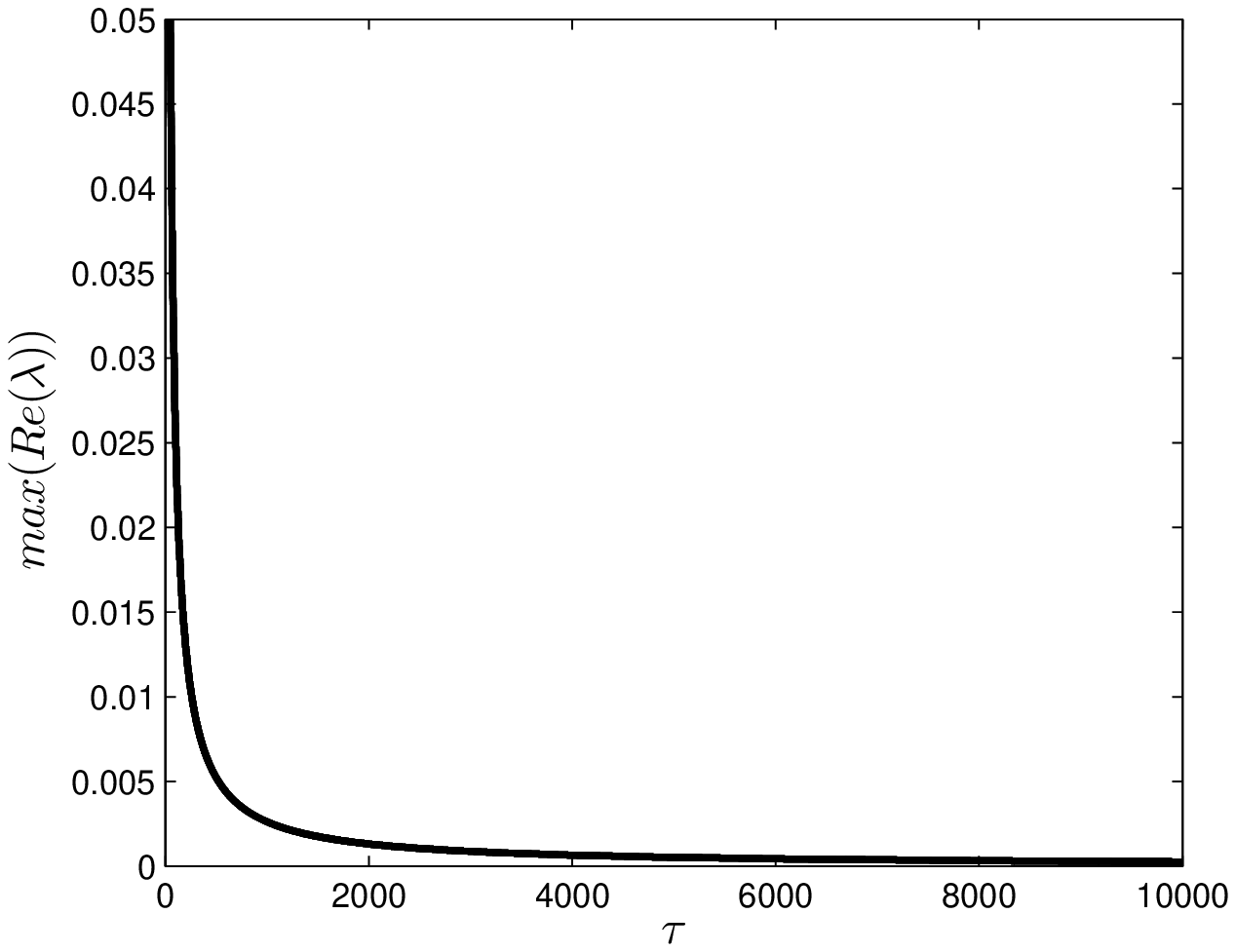}
  \end{subfigure}
  \caption{Real part of the rightmost root for the characteristic function $P_{\Fix(\S_N)}(\lambda,\tau)$ with $K=2$ and $\phi^+$, for $\mu=\{0.1,~0.2,~0.4,~0.6,~0.8\}$.} 
\label{fig:fix_point_full-phase_rightmost_root}
\end{figure}
Using the Matlab routines DDE-Biftool~\cite{Engelborghs2002, Engelborghs2001} we observe that the real parts of the other characteristic roots converge  to $0$ as $\tau\to\infty$, see figure~\ref{fig:fix_point_full-phase_rightmost_root_2}. We are only interested  in finite values of time-delay, consequently there is numerical evidence that some roots in $P_{\Fix(\S_N)}$  remain unstable for $\phi^+$ for any finite value of $\mu,\tau\in\R^+$.
\begin{figure}[!htb]
  \centering
\includegraphics[scale=0.4]{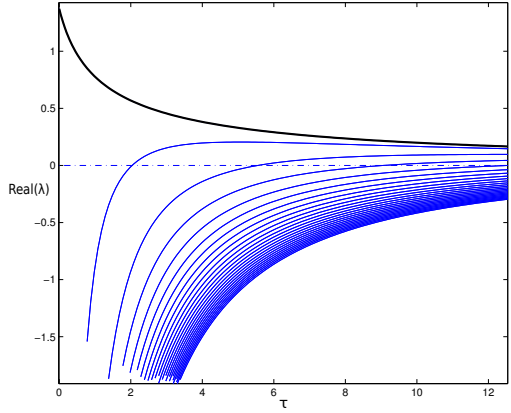}
  \caption{Real part of the rightmost roots for the characteristic function $P_{\Fix(\S_N)}(\lambda,\tau)$ with $K=2$ and $\phi^+$, for $\mu=0.9$ using DDE-Biftool.} 
\label{fig:fix_point_full-phase_rightmost_root_2}
\end{figure}

\item For the equilibrium $\phi^-$ with $K>1$  the characteristic function $P_{\Fix(\S_N)}$ in equation~\eqref{eq:full_phase_P} becomes
\begin{equation*}
P_{\Fix(\S_N)}=\lambda^2+\mu\lambda+K\mu\left(1+\sqrt{1-\dfrac{1}{K^2}}\right)-K\mu\left(1-\sqrt{1-\dfrac{1}{K^2}}\right)e^{-\lambda\tau}=0.
\end{equation*}

From  \eqref{eq:full_phase_charc_fun_tau0_a} we know both roots are stable when $\tau=0$, and from \eqref{eq:full_phase_roots_tau_infty} we also know there are no roots in the right-side of the complex plane when $\tau\to\infty$. 

From $F(\omega)=0$ in~\eqref{eq:full_phase_exist_cond} with $\phi^*=\phi^-$, we obtain
\begin{equation*}
\begin{array}{rcl}
 \omega_{\pm}^2&=&-\dfrac{1}{2}\left(\mu^2-2K\mu\left(1+\sqrt{1-\dfrac{1}{K^2}}\right) \right)\\
&&\pm\dfrac{1}{2}\left[\left(\mu^2-2K\mu\left(1+\sqrt{1-\dfrac{1}{K^2}}\right) \right)^2-16K^2\mu^2\sqrt{1-\dfrac{1}{K^2}}\right]^{1/2}.
\end{array}
\end{equation*}
Then by Lemma \ref{l.conditions_bif_Xj}, using that   $c>0$ by \eqref{e.bc} , we see that  $\omega\in\R^+$ if and only if the discriminant is positive and
 $b<0$, i.e., 
the first term is positive which implies, 
\begin{equation}
  \label{eq:full_phase_cond_C1}
  \mu<2\left(K+\sqrt{K^2-1}\right).
\end{equation}
 The discriminant is positive if and only if
\marginnote{
\tiny
\begin{align}
  \left(\mu^2-2K\mu\left(1+\sqrt{1-\dfrac{1}{K^2}}\right) \right)^2~~~~&\nonumber
\\-16K^2\mu^2\sqrt{1-\dfrac{1}{K^2}}>0&\nonumber
\\\left(\mu-2\left(K+\sqrt{K^2-1}\right) \right)^2~~~~&\nonumber
\\-16K\sqrt{K^2-1}>0&\nonumber
\\\ \Leftrightarrow \left|\mu-2\left(K+\sqrt{K^2-1}\right) \right|>4\sqrt{K\sqrt{K^2-1}}&\nonumber
\end{align}
\normalsize
}
\begin{equation}
\left|\mu-2\left(K+\sqrt{K^2-1}\right) \right|>4\sqrt{K\sqrt{K^2-1}},
\end{equation}
and so, by \eqref{eq:full_phase_cond_C1}
\begin{equation}
         \mu<2\left(K+\sqrt{K^2-1}\right)-4\sqrt{K\sqrt{K^2-1}} =: \mu_{\max},\label{eq:full_phase_cond_total}
 \end{equation}
is a  necessary condition for the existence of bifurcations in $\Fix(\S_N)$ for the equilibrium $\phi^-$ with $K>1$.
 When $K=1$ this condition becomes condition~\eqref{eq:full_phase_other_cond}.
\begin{remark}
We know that   $\phi^-$ is spectrally stable at $\tau=0$  for $K\geq1$, see \eqref{eq:full_phase_charc_fun_tau0_a}; using the time-delay as parameter bifurcations can occur for time delays $\tau$ satisfying \eqref{eq:full_phase_tau_general} provided condition~\eqref{eq:full_phase_cond_total} holds. If condition~\eqref{eq:full_phase_cond_total} does not hold then the equilibrium remains stable for all $\tau$. In that sense $\mu_{\max}$ sets the lower limit to $\mu$ for the  stability of the equilibrium  in $\Fix(\S_N)$  for all time delays
$\tau$ in this case.
\end{remark}

As an example we choose
$K=1.05$. Condition~\eqref{eq:full_phase_cond_total} becomes
$0<\mu<\mu_{\max}=0.4211$. With these parameters figure~\ref{fig:full_phase_Pfix_phi2_bif_curves_a} shows the curves $\tau_{\pm}(\mu)$ ($\tau_+$ as a  solid line and $\tau_-$ as a dashed line) 
determined by ~\eqref{eq:full_phase_w},~\eqref{eq:full_phase_sin_cos_Pfix},
and~\eqref{eq:full_phase_tau_general}, at the equilibrium $\phi^*=\phi^-$,  for different values of $n$ (a lobe $\tau_+(n)$ and~$\tau_-(n)$ for each $n$) considering only positive values of $\tau$; we already know from~\eqref{eq:full_phase_AC+BD_b} that $\delta_+(\omega_+)>0$ and $\delta_-(\omega_-)<0$, thus the  shadowed area indicates the region where the equilibrium is
stable in $\Fix(\S_N)$.
These results are tested using DDE-Biftool~\cite{Engelborghs2002, Engelborghs2001}; with $\mu=0.3$ and $K=1.05$ we found the time delays $\tau_1$ to $\tau_5$ corresponding to critical time delays leading to Hopf bifurcations ($\tau_1=6.34$, $\tau_2=11$, $\tau_3=15.41$, $\tau_4=23.51$, and $\tau_5=24.48$) which match  those found using the $S_n$ map in figure~\ref{fig:full_phase_Pfix_phi2_bif_curves_a}; In figure~\ref{fig:full_phase_Pfix_phi2_DDEBiftool} the real part of the rightmost root  is shown as a  black curve,  the critical time delays $\tau_1$ to $\tau_5$ are also shown, each peak is related to the corresponding lobe in figure~\ref{fig:full_phase_Pfix_phi2_bif_curves_a}; the numerics confirms that at $\tau_1$, $\tau_3$ and $\tau_5$ the root crosses from the left to the right of the imaginary axis switching stability from stable to unstable, and at $\tau_2$, $\tau_4$ the roots come back to the left-hand side of the complex plane, switching stability from unstable to stable again; these time delays are the same as shown in figure~\ref{fig:full_phase_Pfix_phi2_bif_curves_a}, clearly for $\tau>\tau_5$ the equilibrium becomes unstable in $\Fix(\S_N)$. 
Thus, for the given parameters the equilibrium is stable  in 
$\Fix(\S_N)$  within the interval $(0,\tau_1)\cup(\tau_2,\tau_3)\cup(\tau_4,\tau_5)$, when $\mu<\mu_{\max}$.
\begin{figure}[!htb]
  \centering
 \includegraphics[scale=0.5]{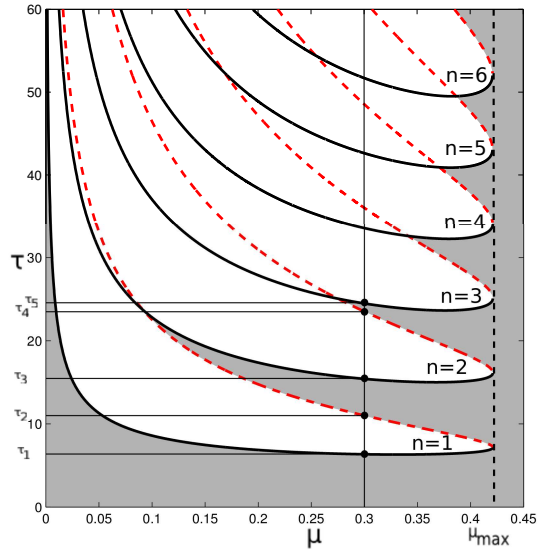}
  \caption{Symmetry-preserving bifurcation curves for the equilibrium $\phi^-$, with $K=1.05$. Within the shadowed regions there are no roots with positive real part.}
  \label{fig:full_phase_Pfix_phi2_bif_curves_a}
\end{figure}
\begin{figure}[!htb]
  \centering
  \includegraphics[scale=0.5]{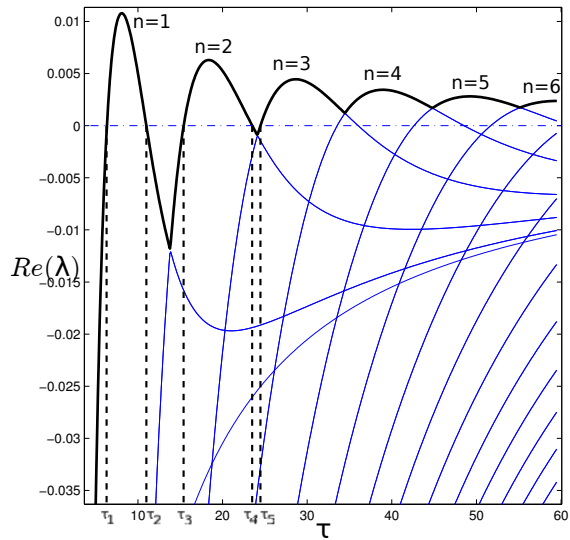}
  \caption{Real part of the rightmost root of $P_{\Fix(\S_N)}$ for $\phi^-$, $\mu=0.3$, and $K=1.05$, using DDE-Biftool.}
\label{fig:full_phase_Pfix_phi2_DDEBiftool}
\end{figure}
\end{itemize}


\subsection{Symmetry-breaking bifurcations}
\label{subsec:Bif_Pj_full_phase_model}

\subsubsection{Roots of the characteristic function $P_U(\lambda,\tau)$ at $\tau=0$ and as $\tau\to\infty$}
Remembering that $K\geq1$, $N\in\N>1$ and $\mu\in\R^+$, we have that the characteristic function $P_U(\lambda,\tau)$ in \eqref{eq:full_phase_P} when $\tau=0$ becomes
\begin{equation*}
  P_U(\lambda,0)=\lambda^2+\mu\lambda+K\mu(1-\cos(2\phi^*))+\frac{K\mu}{N-1}(1+\cos(2\phi^*))=0,
\end{equation*}
where $\phi^*=\phi^\pm$,
which has the two roots
\begin{equation*}
\lambda_{\pm}=-\frac{\mu}{2}\pm\frac{1}{2}\left(\mu^2-4\left[K\mu(1-\cos(2\phi^*))+\frac{K\mu}{N-1}(1+\cos(2\phi^*))\right] \right)^{1/2}.
\end{equation*}
since $|\cos(2\phi^*)|<1$, see \eqref{eq:full_phase_phi**},   the discriminant is always smaller than $\mu^2$, consequently $\Re(\lambda_{\pm})<0$.

When $\tau\to\infty$, assuming $\Re(\lambda)>0$ in \eqref{eq:full_phase_P}, we obtain,
\begin{equation*}
  \lim_{\tau\to\infty}\lambda_{\pm}=-\frac{\mu}{2}\pm\frac{1}{2}\left(\mu^2-4K\mu\left(1-\cos(2\phi^*) \right) \right)^{1/2}.
\end{equation*}
Here again  the discriminant is always smaller than $\mu^2$, thus $\Re(\lambda_{\pm})<0$, which contradicts the assumption $\Re(\lambda)>0$, therefore the roots of $P_U(\lambda,\tau)$ are not in the right-hand side of the complex plane as $\tau\to\infty$.
These results are valid for both equilibria $\phi^\pm$.

\subsubsection{Conditions for the existence of symmetry-breaking bifurcations}
For the characteristic function $P_U(\lambda,\tau)$ from \eqref{eq:full_phase_P}, following \eqref{eq:Sn_map_P} we have,  with $\phi^*=\phi^\pm$,
\begin{equation}\label{e.RSSymBrFullPhase}
  \begin{array}{rcl}
    R(\lambda)&=&\lambda^2+\mu\lambda+K\mu(1-\cos(2\phi^*)) \\
    S(\lambda)&=&\dfrac{K\mu}{N-1}(1+\cos(2\phi^*)),
  \end{array}
\end{equation}
and substituting $\lambda=\i\omega$ we obtain
\begin{equation}
  \label{eq:full_phase_f_Xj_R_S_a}
  \begin{array}{rcl}
    R(\i\omega)&=&-\omega^2+K\mu(1-\cos(2\phi^*))+\i\mu\omega \\
    S(\i\omega)&=&\dfrac{K\mu}{N-1}(1+\cos(2\phi^*)).
  \end{array}
\end{equation}
Then the polynomial $F(\omega)$ from \eqref{eq:Sn_map_F} becomes
\begin{equation}
  \label{eq:full_phase_f_Xj}
  \begin{array}{rcl}
   F(\omega)&=&\omega^4+\left(\mu^2-2K\mu\left(1-\cos(2\phi^*) \right) \right)\omega^2\\
&&+(K\mu)^2\left(1-\cos(2\phi^*)\right)^2-\left(\dfrac{K\mu}{N-1}\right)^2\left(1+\cos(2\phi^*) \right)^2,
  \end{array}
\end{equation}
the roots of which are
\begin{equation}
  \label{eq:full_phase_w_Xj}
  \begin{array}{rcl}
 \omega_{\pm}^2&=&-\dfrac{1}{2}\left(\mu^2-2K\mu\left(1-\cos(2\phi^*) \right) \right)\\
&&\pm\dfrac{1}{2}\Bigg[\left(\mu^2-2K\mu\left(1-\cos(2\phi^*) \right) \right)^2\\
&&\left.-4\left\{(K\mu)^2\left(1-\cos(2\phi^*) \right)^2-\left(\dfrac{K\mu}{N-1}\right)^2\left(1+\cos(2\phi^*) \right)^2 \right\} \right]^{1/2}.
\end{array}
\end{equation}
For sake of simplicity we write
\begin{equation}
\label{eq:full_phase_w_Xj_simple}
\omega_{\pm}=\sqrt{-\frac{b}{2}\pm\frac{1}{2}\sqrt{b^2-4c}},  
\end{equation}
where
\begin{equation}
  \label{eq:full_phase_w_Xj_bc}
  \begin{array}{rcl}
    b&=&\mu^2-2K\mu(1-\cos(2\phi^*))\\
    c&=&(K\mu)^2(1-\cos(2\phi^*))^2-\left(\dfrac{K\mu}{N-1}\right)^2(1+\cos(2\phi^*))^2.
  \end{array}
\end{equation}
The first necessary condition for the existence of symmetry breaking bifurcations   is \eqref{eq:full_phase_Pj_general_exist_cond_a}. 
Substituting $b$ and $c$ into  \eqref{eq:full_phase_Pj_general_exist_cond_a}
 we obtain
\marginnote{
\tiny
\begin{align*}
\left(\mu^2-2K\mu\left(1-\cos(2x_1^*) \right) \right)^2~~~~&
\\-4\left\{(K\mu)^2\left(1-\cos(2x_1^*) \right)^2\right.~~~~&
\\\left.-\left(\dfrac{K\mu}{N-1}\right)^2\left(1+\cos(2x_1^*) \right)^2 \right\}\geq0&
\\\left(\mu-2K\left(1-\cos(2x_1^*) \right) \right)^2~~~~&
\\-4K^2\left(1-\cos(2x_1^*) \right)^2~~~~&
\\+\left(\dfrac{2K}{N-1}\right)^2\left(1+\cos(2x_1^*) \right)^2\geq0&
\\\mu^2-4\mu K(1-\cos(2x_1^*))~~~~&
\\+\left(\dfrac{2K}{N-1}\right)^2(1+\cos(2x_1^*))^2\geq0&
\end{align*}
\normalsize
}
\begin{equation}\label{e.discriminant}
 \begin{array}{rcl}
\mu^2-4\mu K(1-\cos(2\phi^*))+\left(\dfrac{2K}{N-1}\right)^2(1+\cos(2\phi^*))^2&\geq&0.
\end{array}
\end{equation}
This condition is always true if $c \leq 0$ or 
\begin{equation}
(N-1)^2\left(1-\cos(2\phi^*)\right)^2-(1+\cos(2\phi^*))^2 \leq 0,
\end{equation}
which  is equivalent to
\begin{equation}
\label{eq:full_phase_condition_c}
\dfrac{N-2}{N}\leq \cos(2\phi^*).
\end{equation}
\marginnote{
\tiny
\begin{align*}
(K\mu)^2\left(1-\cos(2x_1^*) \right)^2~~~~&
\\-\left(\dfrac{K\mu}{N-1}\right)^2\left(1+\cos(2x_1^*) \right)^2>0&
\\(N-1)^2\left(1-\cos(2x_1^*)\right)^2~~~~&
\\-(1+\cos(2x_1^*))^2>0&
\\(N-1)\left(1-\cos(2x_1^*)\right)-(1+\cos(2x_1^*))>0&
\\\dfrac{N-2}{N}>\cos(2x_1^*)&
\end{align*}
\normalsize
}
If \eqref{eq:full_phase_condition_c} does not hold, then, 
by calculating the real roots of the left hand side of 
\eqref{e.discriminant}, it is possible to find the boundaries in which \eqref{e.discriminant} holds true; 
these are the two curves depending on $K$ and $N$,
\begin{equation}
  \label{eq:full_phase_Pj_general_exist_cond_f}
\begin{array}{l}
  \mu_{\pm}(K;N)=
\\2K(1-\cos(2\phi^*))\pm2K\left[\left(1-\cos(2\phi^*)\right)^2-\left(\dfrac{1}{N-1}\right)^2(1+\cos(2\phi^*))^2\right]^{1/2}.
\end{array}
\end{equation}
Note that the discriminant is always smaller than the square of the first term, and non-negative  for $c\geq 0$. In this case
 $\mu_{\pm}\in\R^+$ and the set $M$ of all values $\mu$ satisfying condition~\eqref{eq:full_phase_Pj_general_exist_cond_a} is
\begin{equation}
  \label{eq:full_phase_Pj_M} 
M= ( 0,\mu_{-}]\cup[\mu_{+},\infty).
\end{equation}
Additional necessary conditions for the existence of Hopf bifurcations 
are given in Lemma
\ref{l.conditions_bif_Xj}.
The condition $b\geq0$ is equivalent to 
\begin{equation}
  \label{eq:full_phase_Pj_general_exist_cond_b}
  \mu\geq2K(1-\cos(2\phi^*)),
\end{equation}
Now we will start the analysis of the conditions for the existence of symmetry breaking Hopf bifurcations  considering three cases:
Now we will start the analysis of the conditions for the existence of bifurcations in $X_j$ considering three cases:
\begin{itemize}
\item When $K=1$ then $\cos(2\phi^*)=0$, see
 \eqref{eq:full_phase_phi**}, and  $\phi^*=\phi^+=\phi^-$. For this case the curves $\mu_{\pm}$ from \eqref{eq:full_phase_Pj_general_exist_cond_f} become
\marginnote{
\tiny
\begin{align*}
\mu_\pm &= 2 \pm 2\left( 1 - \frac{1}{(N-1)^2}\right)^{1/2}\\
& = 
2 \pm \frac{2}{N-1} \sqrt{(N-1)^2-1 }\\
& = 
2 \pm \frac{2}{N-1} \sqrt{N^2-2N } 
\end{align*}
\normalsize
}
\begin{equation}
  \label{eq:full_phase_cond_K=1_a}
  \mu_{\pm}(N)=2\pm\frac{2}{N-1}\sqrt{N(N-2)},
\end{equation}
clearly, for $N\in\N>1$ we have 
\begin{equation}
  \label{eq:full_phase_cond_K=1_b}
0<\mu_{-}(N)\leq2\leq\mu_{+}(N).
\end{equation}
From \eqref{eq:full_phase_condition_c} we see that $c\geq0$ 
is always true for this case. 
From Lemma \ref{l.conditions_bif_Xj}, case 2b), we know
 that $b$ from \eqref{eq:full_phase_w_Xj_bc} has to satisfy $b\leq 0$ for real solutions  $\omega$ of \eqref{eq:full_phase_f_Xj} to exist
which becomes
\begin{equation}
  \label{eq:full_phase_cond_k=1_d}
  \mu\leq 2.
\end{equation}
Thus from \eqref{eq:full_phase_cond_K=1_b},~\eqref{eq:full_phase_cond_k=1_d} and \eqref{eq:full_phase_Pj_M} we see that symmetry breaking bifurcations  can appear if and only if
\begin{equation}
  \label{eq:full_phase_cond_k=1_e}
  \mu\in(0,\mu_{-}].
\end{equation}
\item  For the equilibrium $\phi^+$, with $K>1$ we have $\cos(2\phi^+)=\frac{1}{K}\sqrt{K^2-1}$. We know  that this equilibrium
is unstable in $\Fix(\S_N)$, see sections~\ref{subsub:Condition_existence_bif_in_Pfix} and~\ref{subsubsec:roots_Pfix_tau0_tau_infty}.  
\marginnote{
\tiny
\begin{align*}
\dfrac{N-2}{N}=&~\dfrac{1}{\mathcal{K}_N}\sqrt{\mathcal{K}_N^2-1}
\\\left(\dfrac{N-2}{N}\right)^2\mathcal{K}_N^2=&~\mathcal{K}_N^2-1
\\1=&~\mathcal{K}^2_N\left(1-\left(\dfrac{N-2}{N}\right)^2 \right)
\\1=&~\mathcal{K}_N\left(N^2 - (N-2)^2 \right)^{1/2}/N
\\1=&~\mathcal{K}_N\left(4N-4 \right)^{1/2}/N
\\\mathcal{K}_N=&~\dfrac{N}{2}\sqrt{\dfrac{1}{N-1}}
\end{align*}
\normalsize
}
The curves for $\mu_{\pm}$ in~\eqref{eq:full_phase_Pj_general_exist_cond_f} for this case become
\begin{equation}
  \label{eq:full_phase_Pj_phi1_K>1_a}
 \begin{array}{l}
  \mu_{\pm}(K;N)=
\\2\left(K-\sqrt{K^2-1}\right)~\pm 2\left[\left(K-\sqrt{K^2-1}\right)^2-\dfrac{1}{(N-1)^2}\left(K+\sqrt{K^2-1}\right)^2 \right]^{1/2}.
\end{array}
\end{equation}
Moreover $b=0$ (with $b$ from \eqref{eq:full_phase_w_Xj_bc}) is satisfied if and only if
\begin{equation}
  \label{eq:full_phase_Xj_phi1_k>1_b=0}
  \mu=\mu_{b}(K):=2\left(K-\sqrt{K^2-1}\right),
\end{equation}
and the curve $c=0$ is equivalent to
\begin{equation}
  \label{eq:full_phase_Xj_phi1_k>1_c=0}
  K={K}_N:=\dfrac{N}{2}\sqrt{\dfrac{1}{N-1}}.
\end{equation}
\marginnote{
\tiny
\begin{align*}
\reversemarginpar
&  \sin(\omega_{\pm}\tau)=\dfrac{R_IS_R-S_IR_R}{|S|^2}
\\&~~~                      =\dfrac{R_IS_R}{|S|^2}
\\&~~~                      =\dfrac{\omega_{\pm}(N-1)}{K(1+\cos(2x_1^*))}
\\&\cos(\omega_{\pm}\tau)=-\dfrac{S_IR_I+S_RR_R}{|S|^2}
\\&~~~                      =-\dfrac{S_RR_R}{|S|^2}
\\&                      =\dfrac{\left(\omega_{\pm}^2-K\mu\left(1-\cos(2x_1^*) \right) \right)(N-1)}{K\mu\left(1+\cos(2x_1^*)\right) }.
\\  A&=\Re\left(e^{-i\omega\tau}\left(i\omega S(i\omega,\tau)\right.\right.
\\   &~~~\left.\left.-S'_{\tau}(i\omega,\tau) \right)-R'_{\tau}(i\omega,\tau) \right)
\\   &=\Re\left(e^{-i\omega\tau}i\omega S(i\omega,\tau) \right)
\\   &=\Re\left(e^{-i\omega\tau}i\omega\dfrac{K\mu}{N-1}\left(1+\cos(2x_1^*)\right)\right)
\\   &=\Re\left((\cos(\omega\tau)-i\sin(\omega\tau) )\right.
\\   &~~~\left.\times\left(i\omega\dfrac{K\mu}{N-1}\left(1+\cos(2x_1^*)\right)\right) \right)
\\   &=\dfrac{\omega_{\pm} K\mu}{N-1}\left(1+\cos(2x_1^*) \right)\sin(\omega_{\pm}\tau)
\\   &=\dfrac{\omega_{\pm} K\mu}{N-1}\left(1+\cos(2x_1^*) \right)\dfrac{\omega_{\pm}(N-1)}{K(1+\cos(2x_1^*))}
\\   &=\omega_{\pm}^2\mu
\\  B&=\Im\left(e^{-i\omega\tau}\left(i\omega S(i\omega,\tau)-S'_{\tau}(i\omega,\tau) \right)\right.
\\   &~~~\left.-R'_{\tau}(i\omega,\tau) \right)
\\   &=\Im\left(e^{-i\omega\tau}\left(i\omega S(i\omega,\tau) \right) \right)
\\   &=\Im\left((\cos(\omega\tau)-i\sin(\omega\tau) )\right.
\\   &~~~\left.\times\left(i\omega\dfrac{K\mu}{N-1}\left(1+\cos(2x_1^*)\right)\right) \right)
\\   &=\dfrac{\omega_{\pm} K\mu}{N-1}\left(1+\cos(2x_1^*) \right)\cos(\omega_{\pm}\tau)
\\   &=\dfrac{\omega_{\pm} K\mu}{N-1}\left(1+\cos(2x_1^*) \right)
\\   &~~~\times\dfrac{\left(\omega_{\pm}^2-K\mu\left(1-\cos(2x_1^*) \right) \right)(N-1)}{K\mu\left(1+\cos(2x_1^*)\right) }
\\   &=\omega_{\pm}\left(\omega_{\pm}^2-K\mu\left(1-\cos(2x_1^*) \right) \right)
\end{align*}
\normalsize
}
In figure~\ref{fig:full_phase_Pj_phi1_exist_cond_curves}, the curves
from
\eqref{eq:full_phase_Pj_phi1_K>1_a},~\eqref{eq:full_phase_Xj_phi1_k>1_b=0}
and~\eqref{eq:full_phase_Xj_phi1_k>1_c=0} are shown for various values
of $N$. The curves $\mu_{+}(K;N)$ are shown in dotted lines and the
curves $\mu_{-}(K;N)$ in dashed lines; the curve $\mu_{b}(K)$ is the
solid black line. Here the conditions $b>0$ and $b<0$ correspond to
the regions above and below the curve $\mu_{b}$ respectively; the
conditions $c<0$ and $c>0$ are identified with the right and  left sides of each vertical line $ {K}_N$. Let us take as an example $N=2$, from \eqref{eq:full_phase_Xj_phi1_k>1_c=0} we see that $ {K}_2=1$, but we know that $K>1$  by assumption, thus bifurcation can occur at the right-hand side of the vertical line  $K= {K}_2$,  where $c<0$. From the additional conditions for the existence of bifurcations given in Lemma \ref{l.conditions_bif_Xj}, we see that for $c<0$ bifurcations occur only for $\omega_+$.   
\begin{figure}[Ihtb]
  \centering
  \includegraphics[scale=0.5]{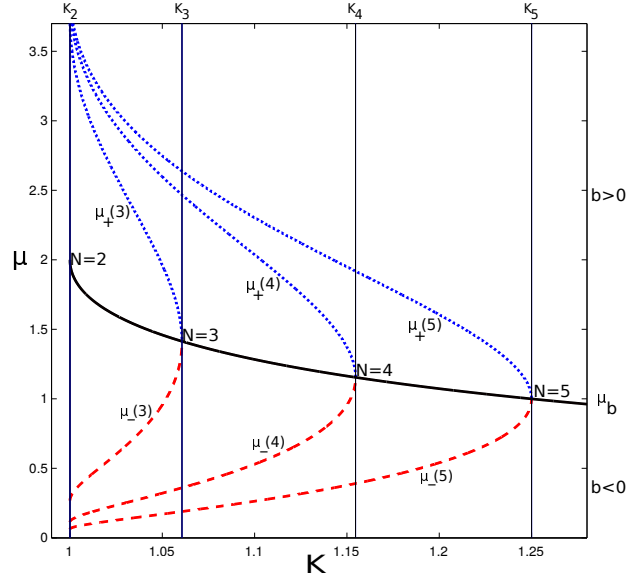}
 \caption{Curves showing the conditions for existence of symmetry-breaking bifurcations  for $\phi^+$ and $K>1$.}
  \label{fig:full_phase_Pj_phi1_exist_cond_curves}
\end{figure}
\item 
\item For the equilibrium $\phi^-$ with $K>1$  the curves $\mu_{\pm}$ from  \eqref{eq:full_phase_Pj_general_exist_cond_f} become
\begin{equation}
  \label{eq:full_phase_X_j_phi2_K>1}
  \begin{array}{l}
  \mu_{\pm}=
\\2\left(K+\sqrt{K^2-1}\right)\pm2\left[\left(K+\sqrt{K^2-1}\right)^2-\dfrac{1}{(N-1)^2}\left(K-\sqrt{K^2-1}\right)^2 \right]^{1/2},
  \end{array}
\end{equation}
and the condition $c>0$  from \eqref{eq:full_phase_condition_c} becomes
\begin{equation}
  \frac{N-2}{N}>-\frac{1}{K}\sqrt{K^2-1}
\end{equation}
which is always true since $K>1$ and $N\in\N>1$. From the conditions for the existence of symmetry breaking bifurcations given in  Lemma \ref{l.conditions_bif_Xj} we see that for $c>0$ we need $b<0$ in order for bifurcations to occur, i.e.,
\begin{equation}
  \label{eq:full_phase_X_j_phi2_K>1_b}
  \mu<2(K+\sqrt{K^2-1}),
\end{equation}
therefore from \eqref{eq:full_phase_X_j_phi2_K>1} and ~\eqref{eq:full_phase_X_j_phi2_K>1_b} we see that symmetry breaking  bifurcations can occur  at the equilibrium $\phi^-$  for both $\omega_{\pm}$ with $K>1$ if
\begin{equation}
  \mu\in(0,\mu_{-}).
\end{equation}
\end{itemize}
The analysis of the roots of the polynomial $F(\omega)$ from \eqref{eq:full_phase_f_Xj} in the above paragraphs gives us necessary conditions for the existence of roots  of $P_U(\i \omega,\tau)$, where $\omega\in\R^+$. However it is necessary to also impose  the conditions on $\sin(\omega\tau)$ and $\cos(\omega\tau)$ given in section~\ref{subsec:Sn_map}, to determine at which time delays the linearization at the equilibrium has imaginary eigenvalues. From \eqref{eq:Sn_map_sin_cos} and~\eqref{eq:full_phase_f_Xj_R_S_a} we have 
\begin{equation}
  \label{eq:full_phase_Xj_sin_cos}
  \begin{array}{rcl}
    \sin(\omega_{\pm}\tau)&=&\dfrac{\omega_{\pm}(N-1)}{K(1+\cos(2\phi^*))}\\\\
    \cos(\omega_{\pm}\tau)&=&\dfrac{\left(\omega_{\pm}^2-K\mu\left(1-\cos(2\phi^*) \right) \right)(N-1)}{K\mu\left(1+\cos(2\phi^*)\right) }.
  \end{array}
\end{equation}
\marginnote{
\tiny
\begin{align*}
 C&=\Re\left(R'_{\lambda}(i\omega,\tau)\right.
\\   &\left.~~~+e^{-i\omega\tau}\left(S'_{\lambda}(i\omega,\tau)-\tau S(i\omega,\tau) \right)  \right)
\\   &=\Re\left(i2\omega+\mu\right.
\\   &~~~\left.+e^{-i\omega\tau}\left(-\tau \dfrac{K\mu}{N-1}(1+\cos(2x^*_1 )) \right) \right)
\\   &=\Re\left(i2\omega+\mu -(\cos(\omega\tau)-i\sin(\omega\tau))\right.
\\   &~~~\left.\times\dfrac{K\mu\tau}{N-1}(1+\cos(2x^*_1 ))\right)
\\   &=\mu-\dfrac{K\tau\mu}{N-1}\left(1+\cos(2x_1^*) \right)\cos(\omega\tau)
\\   &=\mu-\dfrac{K\tau\mu}{N-1}\left(1+\cos(2x_1^*) \right)
\\   &~~~\times\dfrac{\left(\omega_{\pm}^2-K\mu\left(1-\cos(2x_1^*) \right) \right)(N-1)}{K\mu\left(1+\cos(2x_1^*)\right) }
\\   &=\mu-\tau\left(\omega_{\pm}^2-K\mu\left(1-\cos(2x_1^*) \right) \right)
\\D&=\Im\left(R'_{\lambda}(i\omega,\tau)\right.
\\   &~~~\left.+e^{-i\omega\tau}\left(S'_{\lambda}(i\omega,\tau)-\tau S(i\omega,\tau) \right)  \right)
\\   &=\Im\left(i2\omega+\mu\right.
\\   &~~~\left.+e^{-i\omega\tau}\left(-\tau \dfrac{K\mu}{N-1}(1+\cos(2x_1^* )) \right) \right)
\\   &=\Im\left(i2\omega+\mu -(\cos(\omega\tau)-i\sin(\omega\tau))\right.
\\   &~~~\left.\times\dfrac{K\mu\tau}{N-1}(1+\cos(2x_1^* ))\right)
\\   &=2\omega_{\pm}+\dfrac{K\tau\mu}{N-1}\left(1+\cos(2x_1^*)\right)\sin(\omega\tau)
\\   &=\omega_{\pm}(2+\tau\mu).
\end{align*}
\normalsize
}
Note that the denominator in those terms does not vanish for all
 $K,\mu\in\R^+$, $K\geq1$, $N\in\N>1$, since $\cos(2\phi^*)$ from \eqref{eq:full_phase_phi**} satisfies  $|\cos(2\phi^*)|<1$. The frequency $\omega_{\pm}$ is computed using~\eqref{eq:full_phase_w_Xj}. At this point we can calculate the time delays $\tau\geq0$ associated to $\omega_{\pm}$ using \eqref{eq:2n_tau_map}, which for this case becomes,
\begin{equation}
  \label{eq:full_phase_Xj_tau_eq}
  \tau_{\pm}(\mu,K;n)=\frac{1}{\omega_{\pm}}\left(\arg\left( \omega_{\pm}^2-K\mu\left(1-\cos(2\phi^*)\right),  \omega_{\pm}\mu \right)+ 2n\pi \right),~~~n\in\mathbb{N}.
\end{equation}
The last necessary condition for the existence of bifurcation points is the transversality condition $\delta \neq 0$ where $\delta$ is as in ~\eqref{eq:2n_delta}
and \eqref{eq:Sn_delta_abcd}, 
and, due to  \eqref{e.RSSymBrFullPhase} and~\eqref{eq:full_phase_Xj_sin_cos},  $A$, $B$, $C$ and $D$ are as 
in \eqref{eq:full_phase_delta_abcd_general}. As in the case of symmetry-preserving
bifurcations,
we see from   \eqref{eq:full_phase_w_Xj}, \eqref{eq:full_phase_w_Xj_simple} and \eqref{eq:full_phase_w_Xj_bc}  that \eqref{eq:full_phase_AC+BD} holds true again and that therefore the sign of $\delta$ is again given by \eqref{eq:full_phase_AC+BD_b}.
Hence, whenever  $\lambda=\i\omega_-$, $\omega_-\in\R^+$ is a root of $P_U$, then
 it crosses the imaginary axis from the right to the left, whereas whenever $\lambda=\i\omega_+$, $\omega_+\in\R^+$ is a root it   crosses the imaginary axis from the left to the right.
 

\subsubsection{Curves of symmetry-breaking bifurcations}
In the previous section we have analyzed  the conditions for the existence of symmetry-breaking bifurcations in terms of the parameters $\mu,K\in\R^+$ with $K\geq1$ and $N\in\N>1$ for both equilibria $\phi^\pm$. In this section we will obtain curves of symmetry-breaking bifurcations using the $S_n$ map, see section~\ref{subsec:Sn_map}, and we shall compare these curves  with those obtained for  $\Fix(\S_N)$ in section~\ref{subsec:Bif_Pfix_full_phase_model}.
\marginnote{
 \reversemarginpar
 \tiny
 \begin{align*}
 &AC+DB=
 \\&=(\omega_{\pm}^2\mu)(\mu-\tau(\omega_{\pm}^2-K\mu(1-\cos(2x_1^*))))
 \\&+\omega_{\pm}^2(2+\tau\mu)(\omega_{\pm}^2-K\mu(1-\cos(2x_1^*)))
 \\&=\omega_{\pm}^2\mu^2-\mu\tau\omega_{\pm}^2(\omega_{\pm}^2-K\mu(1-\cos(2x_1^*)))
 \\&+\omega_{\pm}^2(2+\tau\mu)(\omega_{\pm}^2-K\mu(1-\cos(2x_1^*)))
 \\&=\omega_{\pm}^2\mu^2+(\omega_{\pm}^2-K\mu(1-\cos(2x_1^*)))(2\omega_{\pm}^2)
 \\&=\omega_{\pm}^2(\mu^2+2(\omega_{\pm}^2-K\mu(1-\cos(2x_1^*))))
 \\&\text{but}~~b=\mu^2-2K\mu(1-\cos(2x_1^*))
 \\&=\omega_{\pm}^2(b+2\omega_{\pm}^2)
 \\&\text{and}~~2\omega_{\pm}^2=-b\pm\sqrt{b^2-4c}
 \\&=\pm\omega_{\pm}^2\sqrt{b^2-4c}
 \end{align*}
 \normalsize
 }

We will consider three cases:
\begin{itemize}
\item When $K=1$ the results we obtain are valid for the unique equilibrium $\phi^+=\phi^- \in (-\pi,\pi]$. We already know that bifurcations for $K=1$ can occur only for $\mu<2$, see \eqref{eq:full_phase_cond_k=1_d}; thus the frequency $\omega$ in  \eqref{eq:full_phase_w_Xj} becomes 
\begin{equation}
  \label{eq:full_phase_X_j_bif_cur_K=1_omega}
  \omega_{\pm}^2=-\frac{1}{2}\mu(\mu-2)\pm\frac{1}{2}\mu\left[(\mu-2)^2-4\left(1-\frac{1}{(N-1)^2} \right) \right]^{1/2}.
\end{equation}
Here we distinguish two cases:
\begin{itemize}
\item When $N=2$ we have  $\omega_+=\sqrt{\mu(2-\mu)}$ and
  $\omega_-=0$ (but $\omega_-=0$ does not correspond to a root of
  $P_U$, so we   ignore it), and from \eqref{eq:full_phase_AC+BD_b} we know that bifurcations associated to $\omega_+$ cross the imaginary axis from the left to the right.
Then plotting the curves for $\tau_+(\mu;n)$ using~\eqref{eq:full_phase_Xj_tau_eq}, and comparing them with those curves obtained for  $\Fix(\S_N)$, see figure~\ref{fig:mu1_x_tau_full_phase}, we obtain the curves shown in  figure~\ref{fig:full_phase_Pj_phi1_Sn_N2}. The curves for $\tau_+$ corresponding to symmetry
breaking bifurcations  are plotted as solid lines and the curves for $\tau_+$ corresponding to bifurcations in  $\Fix(\S_N)$ are  dashed lines, with the time-delay as bifurcation parameter.
 Note that each curve indicates a new root crossing the imaginary axis from the left to the right.
 Therefore within the shadowed region there are no roots in the
 right-hand side of the complex plane, however due to the   zero
 eigenvalue of the linearization in $\Fix(\S_N)$ we cannot conclude nonlinear stability of the equilibrium in this region.
\begin{figure}[!htb]
  \centering
  \includegraphics[scale=0.4]{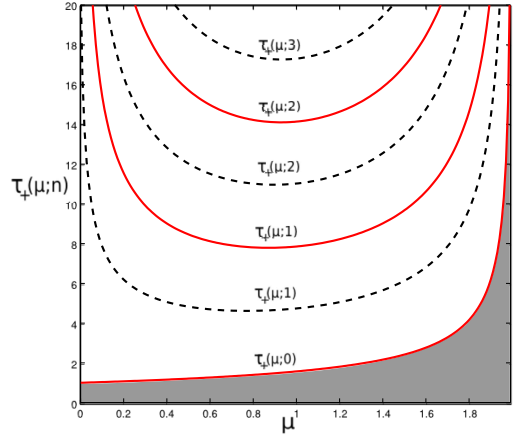}
  \caption{Symmetry-preserving bifurcation curves  in dashed lines and
    symmetry-breaking bifurcation curves in solid line for the 
    equilibrium  $\phi^+=\phi^-$ with $N=2$ and $K=1$; bifurcations
    occur with roots crossing the imaginary axis from the left to the
    right. There are no roots in the right-hand side of the complex
    plane within the shadowed region. However there is a constant zero root
    due to a saddle node bifurcation.}
  \label{fig:full_phase_Pj_phi1_Sn_N2}
\end{figure}
In figure~\ref{fig:full_phase_Pj_phi1_rm_N2} the real part of the
rightmost root for the case $N=2$, $K=1$ and $\mu=\{0.05,~0.5,~1\}$
was computed as a function of the time-delay using the Lambert W
function~\cite{Corless1996, Wang2008}, with both Newton's and Halley's schemes; it can be seen that this root crosses the imaginary axis at a low value of $\tau$ (approx. $\tau=1$), and never comes back, but it approaches zero as $\tau\to\infty$. 
\begin{figure}[!htb]
  \centering
  \begin{subfigure}[t]{2.2in}
  \centering
  \includegraphics[scale=0.4]{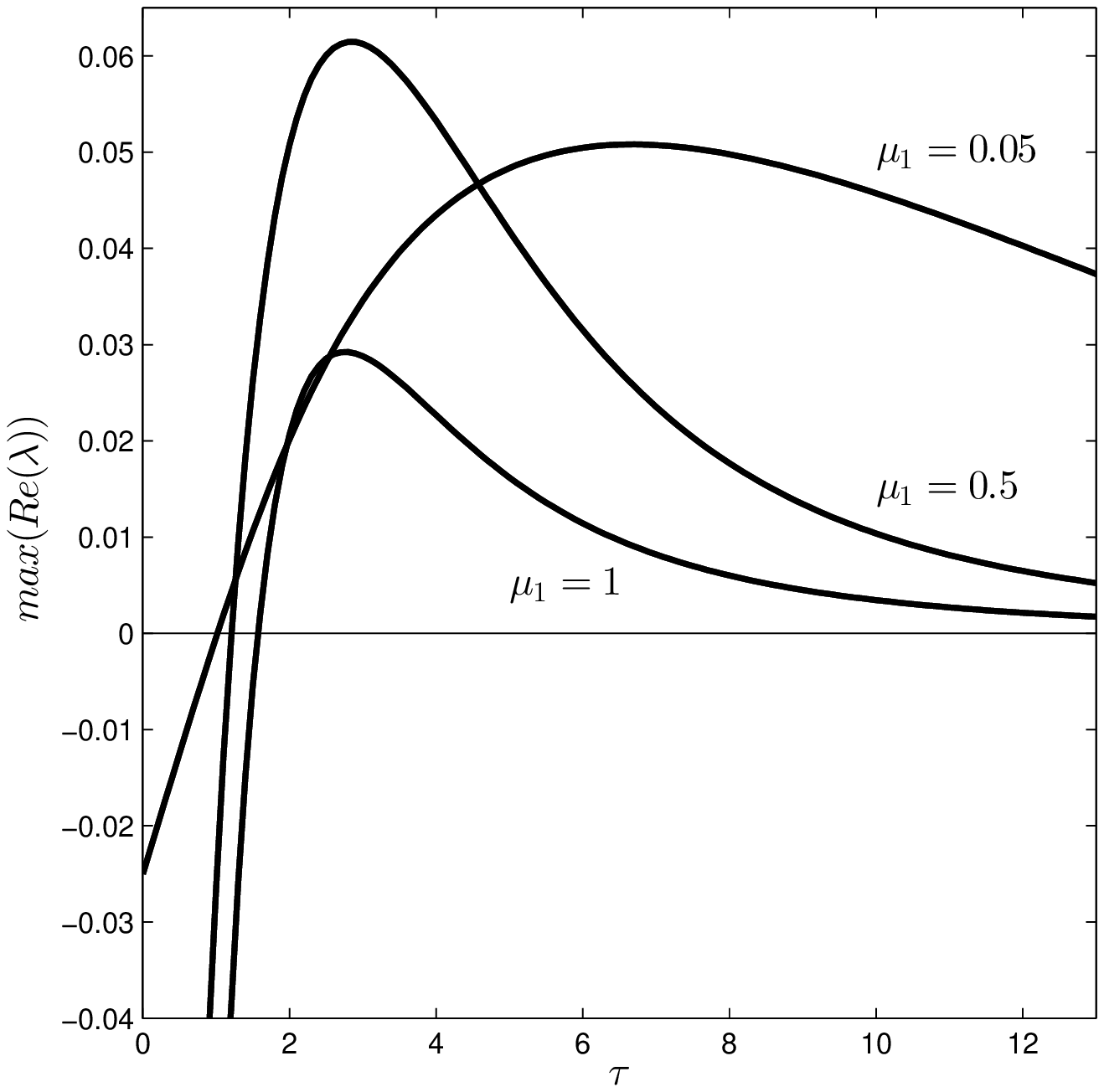}
  \end{subfigure}
  \quad
  \begin{subfigure}[t]{2.2in}
  \centering
  \includegraphics[scale=0.4]{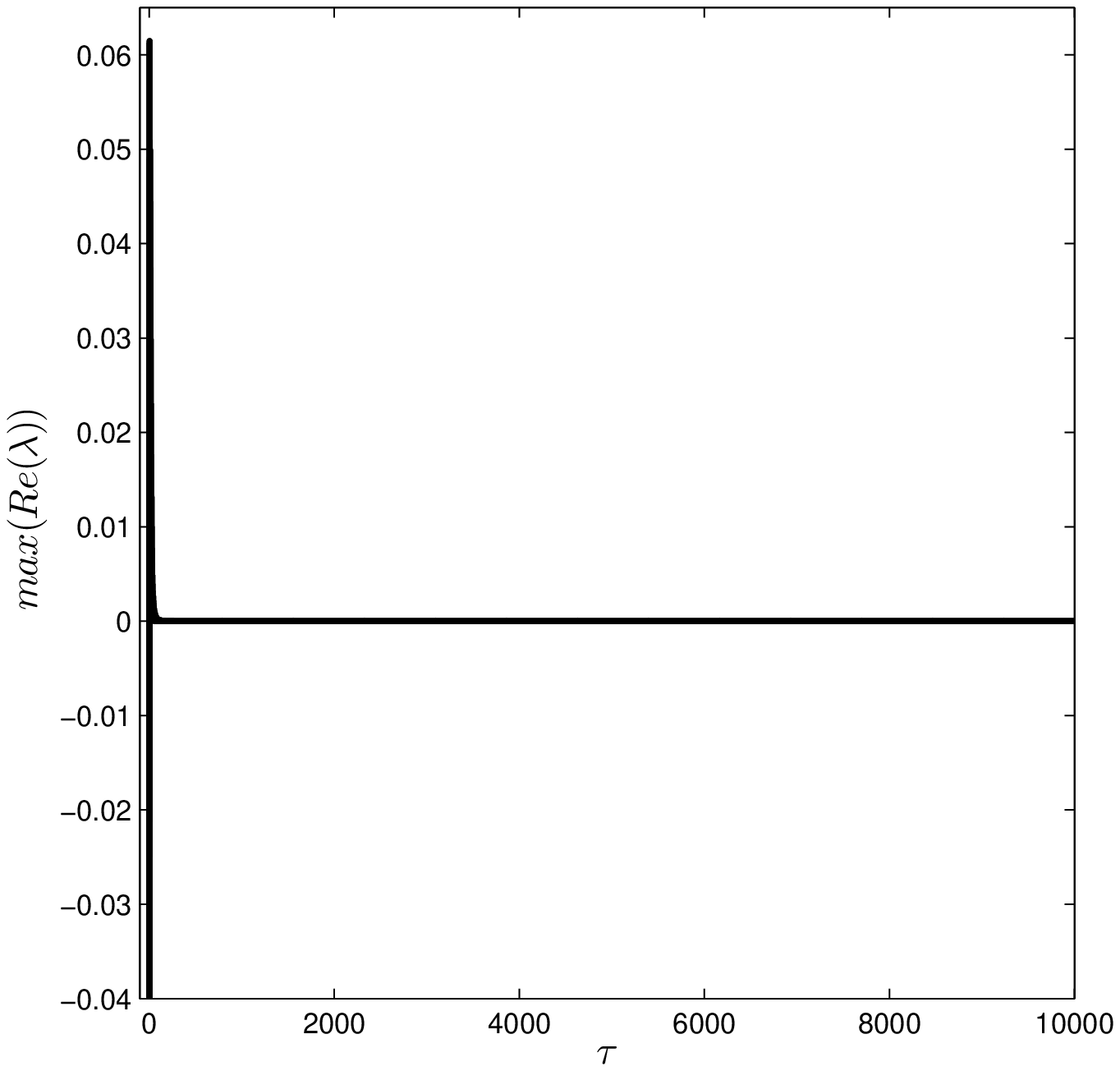}
  \end{subfigure}
  \caption{Real part of the rightmost root for the function $P_U$, for $N=2$, $K=1$ and $\mu=\{0.05,~0.5,~1\}$.}
  \label{fig:full_phase_Pj_phi1_rm_N2}
\end{figure}
\item When $N>2$  we know from  ~\eqref{eq:full_phase_cond_K=1_a}
  and~\eqref{eq:full_phase_cond_k=1_e} that bifurcations can occur for
  $\mu\leq\mu_{1-}(N)<2$, and from
  \eqref{eq:full_phase_X_j_bif_cur_K=1_omega} we see that both
  $\omega_{\pm}\in\R^+$. We also know by looking at the sign of
  $\delta$ in \eqref{eq:full_phase_AC+BD_b} the direction in which these roots cross the imaginary axis as $\tau$ is varied. In   figure~\ref{fig:full_phase_Pj_N3_K=1} the curves of symmetry-breaking bifurcations $\tau_{\pm}$ as a function of $\mu$ from \eqref{eq:full_phase_Xj_tau_eq} with $N=3$ and the curves of symmetry preserving bifurcations   from  \eqref{eq:full_phase_tau} are shown.  The curves of symmetry preserving bifurcations are shown as a solid line. As we  saw in 
section~\ref{subsub:Bifurcation_curves_Fix} these roots cross the imaginary axis from the left to the right. The curves of symmetry-breaking bifurcations  are shown as a dotted line for $\tau_{+}$ and as a dashed line for $\tau_{-}$. The value of $\mu_{-}(N)$ for $N=3$  is also shown, see \eqref{eq:full_phase_cond_K=1_a}, bounding the curves $\tau_{\pm}$. Within the shadowed region there are no roots with positive real part; however there is a zero root of  $P_{\Fix(\S_N)}$.
\begin{figure}[!htb]
  \centering
  \includegraphics[scale=0.45]{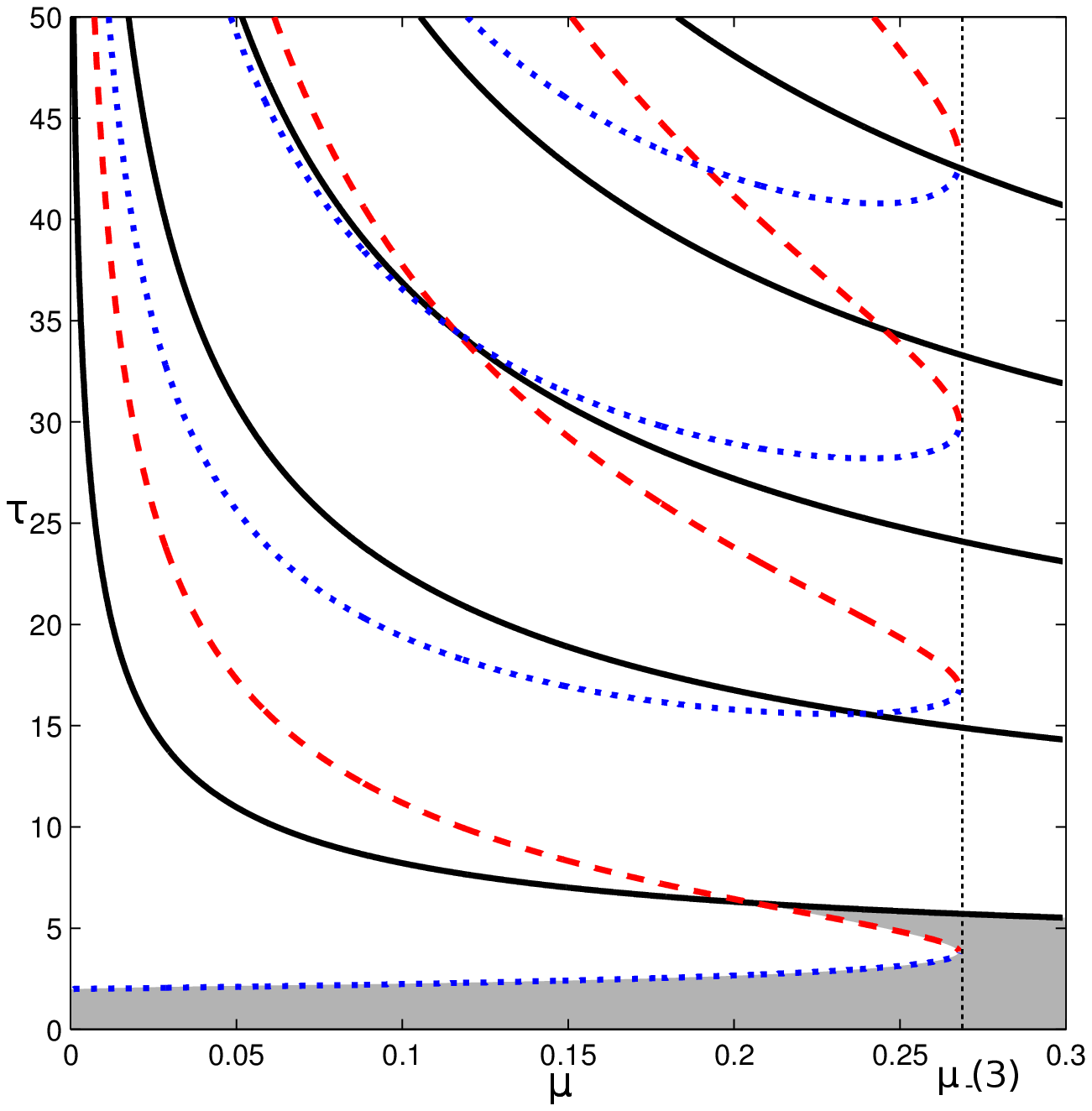}
  \caption{Curves of symmetry-preserving bifurcation as solid line,
    see~\eqref{eq:full_phase_tau}, and 
    curves of symmetry-breaking bifurcations as dotted line for $\mu_{+}$ and as dashed line for
    $\mu_{-} $, see \eqref{eq:full_phase_Pj_phi1_K>1_a}. These curves
    are valid for the equilibrium $\phi^*=\phi^+=\phi^-$, for $N=3$ and $K=1$. Within the shadowed region the equilibrium is spectrally stable.}
  \label{fig:full_phase_Pj_N3_K=1}
\end{figure}
\end{itemize}
\item  If  $K>1$ the equilibrium $\phi^+$ is, as we saw in section~\ref{subsubsec:roots_Pfix_tau0_tau_infty},
unstable  in $\Fix(\S_N)$, nonetheless we can find symmetry-breaking bifurcations of this unstable equilibrium. From Lemma \ref{l.conditions_bif_Xj} we see that
\begin{itemize}
\item If $c\leq0$ then $\omega_{+}\in\R_0^+$.
\item If $c>0$ and $b<0$ then $\omega_{\pm}\in\R^+$.
\end{itemize}
Because of   \eqref{eq:full_phase_w_Xj_bc} the condition $b<0$ implies
\begin{equation*}
  \mu<2(K-\sqrt{K^2-1})<2,
\end{equation*}
and the condition $c>0$ is equivalent to
\begin{equation*}
  K<\ {K}_N,
\end{equation*}
where ${K}_N$ is defined in ~\eqref{eq:full_phase_Xj_phi1_k>1_c=0}. 
For a given $N$ we  see that for small values of $\mu$ and $K$ bifurcations associated to $\omega_{\pm}$ (roots crossing the imaginary axis in both directions) are possible, however for $K\geq {K}_N$ only bifurcations related to $\omega_+$ appear, i.e., roots crossing the imaginary axis from the left to the right. 
\item The equilibrium $\phi^-$ with $K>1$ is  stable in
  $\Fix(\S_N)$ at $\tau=0$, see
  section~\ref{subsubsec:roots_Pfix_tau0_tau_infty}. In
  figure~\ref{fig:full_phase_Pj_N3_Pf_phi2} the curves
  $\tau_{\pm}(\mu)$ for both symmetry preserving and symmetry breaking 
  bifurcations are
  shown; as before the sign of $\delta(\omega_\pm)$ is given
  by~\eqref{eq:full_phase_AC+BD_b}.
\begin{figure}[!htb]
  \centering
  \includegraphics[scale=0.45]{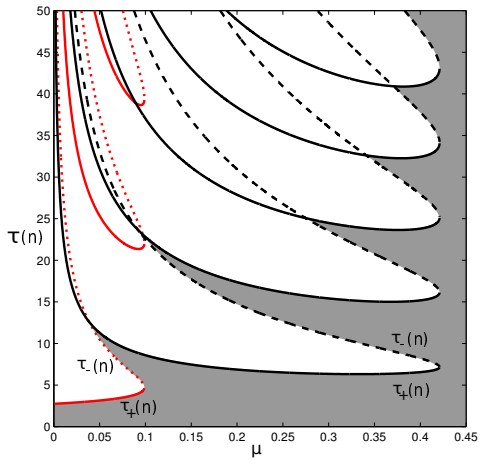}
  \caption{For the equilibrium $\phi^-$ with $N=3$ and $K=1.05$ curves
of symmetry-preserving bifurcations are shown  on the right side  in dashed/solid black lines and curves of  symmetry-breaking bifurcations are shown  on the left side in dotted/solid red lines. Within the shadowed regions the system remains stable.}
  \label{fig:full_phase_Pj_N3_Pf_phi2}
\end{figure}
\end{itemize}


\subsection{Equivariant Hopf bifurcation for $N=2,3$}
\label{subsec:spatio-temp_symm}

Assume~\eqref{eq:RFDE} has a periodic solution $x(t)$ with period $T$. There are two types of symmetry that leave the solution invariant. The first one is the group of \textit{spatial symmetries}
\begin{equation}
  \mathcal K=\{\gamma\in\Gamma, \gamma x(t)=x(t),~\text{for all}~t\},
\end{equation}
which is the isotropy group of each point on the solution. The second is the group of \textit{spatio-temporal symmetries}
\begin{equation}
  \mathcal H=\{\gamma\in\Gamma, \gamma x(t)=x(t+t_0(\gamma)T),~\text{for all}~t\}
\end{equation}
where $t_0(\gamma)\in \R/\Z\cong \S^1$.  
 
The full-phase model \eqref{e.fullPhase}  for two nodes has  $\S_2 = \Z_2$ as symmetry group.
In this case  purely imaginary roots of  
$P_{\Fix(\S_2)}$  lead to symmetry-preserving Hopf bifurcation of fully-synchronized periodic orbits with
spatial symmetry group $\mathcal K=\S_2$. Purely imaginary roots of  
$P_{U}$  lead to symmetry-breaking Hopf bifurcation  of periodic orbits with $\Z_2$ spatio-temporal symmetry
and trivial spatial symmetry,  i.e., the first and second oscillator are half a period out of phase. This follows
from the Equivariant Hopf Theorem,  for details see \cite{Golubitsky1988}.

\begin{figure}[!htb]
  \centering
  \includegraphics[scale=0.45]{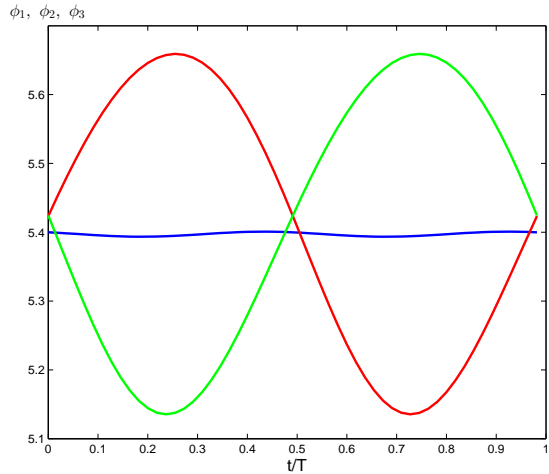}
   \caption{Periodic solutions with spatio-temporal symmetry $\Z_2(\pi_{12})$, for $K=1.05$, $N=3$
    $\mu=0.075$, at $\tau=9.5$, with period $T=24.1895$, bifurcating from the equilibrium $\phi^-=5.3429$ at
    $\tau=7.4898$ with  Hopf frequency $\omega=0.2942$.  Here $ \phi_1$ is  in blue, $ \phi_2$  is in red, and $ \phi_3$ 
is in green, for details see text.}
  \label{fig:symmetry_breaking_psol_phi1_phi2_phi3}
\end{figure}

The full-phase model \eqref{e.fullPhase}  for three nodes has $\S_3 = \DD_3$  as symmetry group  where 
$\DD_m$ is the \textit{dihedral group} of order $2m$ (rotation and
reflections in the plane). In this case  purely imaginary roots of  
$P_{\Fix(\S_3)}$  lead to symmetry-preserving Hopf bifurcation of fully-synchronized periodic orbits with
spatial symmetry group $\mathcal K=\S_3$ as before.
Purely imaginary roots of  
$P_{U}$  lead to symmetry-breaking Hopf bifurcation  of  three families of  periodic orbits (modulo symmetry),
one  with $\mathcal H=\Z_3$  as spatio-temporal symmetry group, corresponding to coordinate shifts $i\to i+1 \mod 3$
and trivial spatial symmetry group,   one with  with $\mathcal K=\Z_2(\pi_{1,2})$-spatial symmetry and  one with spatio-temporal symmetry group $\mathcal H= \Z_2(\pi_{1,2})$, for details see \cite{Golubitsky1988}. For the first family of periodic orbits
the oscillators are 1/3 of a period out of phase. In the second family of periodic orbits, the first and second oscillator are in phase and all oscillators have the same period, in the third family, the  second and first oscillator are out of phase by half a period and the third one  oscillates with twice the period of the other ones, see Figure \ref{fig:symmetry_breaking_psol_phi1_phi2_phi3}.
For systems with $\S_N$-symmetry where $N>3$  the spatial and spatio-temporal symmetry groups of
bifurcating periodic orbits in the case of equivariant Hopf bifurcation have been classified in  \cite{Dias2009}.
%


\section{Phase model}
\label{sec:Phase_model}
In this section a bifurcation analysis for the model of a fully connected N-node network of second-order PLLs oscillators using the instantaneous phase is addressed. 
The phase model has been widely used to analyze the dynamics of
PLL networks for  decades. As a very short review we  mention that in~\cite{Piqueira2003} the influence of the individual
gain of the nodes on the synchronous state is explored,
in~\cite{Bueno2009}  it is analyzed how the filter order influences the
admissible number of nodes in order to reach synchronization,
in~\cite{Bueno2010}  a method is proposed to calculate the PLL filter
in order to successfully suppress the double-frequency term,
in~\cite{Piqueira2006b} the influence of the double-frequency term in
a master-slave strategy is addressed, and in~\cite{Piqueira2005a} the
limitation of a master-slave PLL network is analyzed.

Here, as in the classical approach to the PLL network, we neglect the double-frequency term   and use
 the instantaneous phase $\theta(t)$  instead of the full-phase
 $\phi(t)$ to find time-delays which lead to bifurcation. The model for the $i$-th node
from ~\cite{Bueno2009, Piqueira2006} is
\begin{equation}
    \label{eq:N-node_2nd-order_network}
\ddot{\theta}_i(t)+\mu\dot{\theta}_i (t)-\frac{K\mu}{N-1}\sum_{j=1,j\neq i}^N\sin(\theta_j(t-\tau)-\theta_i(t)-\omega_M\tau)=0,
~~i=1,\ldots, N.
  \end{equation}
This model presents $\S_N\times\sS^1$-symmetry; the demonstration for the $\S_N$-symmetry is similar to the full-phase model and will be omitted here. The translational $C\in\sS^1 = \R/(2\pi \Z)$ acts on $\theta\in \T^N$, $\T^N = (\R/2\pi \Z)^N$ as $\theta_j\to(\theta_j+C) \mod 2\pi$, $j=1,\ldots,N$, and it is not difficult to see that if $\theta(t)=[\theta_1(t),\ldots,\theta_N(t)]^T$ is a solution to~\eqref{eq:N-node_2nd-order_network} then any translated solution
$\tilde\theta(t)$  with  $\tilde\theta_j(t) = \theta_j(t)+C$, for $j=1,\ldots, N$, is also a solution which proves that \eqref{eq:N-node_2nd-order_network}
has translation symmetry. 
This system presents relative equilibria.

\begin{definition}[Relative Equilibrium]
A relative equilibrium of a $\Gamma$-equivariant dynamical system $\dot X = F(X)$ is a group orbit that is invariant under the dynamics.
A trajectory $X(t)$ lies on a relative equilibrium if and only if it is an equilibrium in a comoving frame which 
generates a one-parameter family $\gamma_t$, $t\in \R$, of symmetries, see, e.g., ~\cite{Fiedler1996}.
\end{definition}

 Equation \eqref{eq:N-node_2nd-order_network} has $\S_N$-invariant relative equilibria
\begin{equation}\label{e.RE}
   \theta_k(t)=\Omega(\tau)t+\theta^0,~~~k=1,\ldots,N,
\end{equation}
where the rotation frequency $\Omega(\tau)$ satisfies
\begin{equation}
  \label{eq:Omega(tau)_definition}
   \Omega(\tau)=-K\sin((\Omega(\tau)+\omega_M)\tau),
\end{equation}
and $\theta^0$ is an arbitrary constant; although $\Omega$ depends on $\omega_M$, $K$ and $\tau$, we write $\Omega(\tau)$ because we are interested in the time-delay as bifurcation parameter.

The $\sS^1$-symmetry introduced by the simplification of the double-frequency term generates a zero eigenvalue in the characteristic function for the linearized system around any equilibrium or relative equilibrium.
The $\S_N$-symmetry allows us to find a block decomposition of the linearization at an
$\S_N$-invariant  relative equilibrium as we did in section~\ref{subsec:SN_symm_full_phase_model}.
\begin{remark}
 Relative equilibria which are not $\S_N$-invariant might exist if
 $\theta_i(t)=\Omega(\tau)t+\theta^0_i$ for all $i$, with $\theta^0_i\neq\theta^0_j$, but this case is not studied here. 
\end{remark}
In a rotating frame with frequency $\Omega(\tau)$ such that
\begin{equation*}
  \begin{array}{l}
    \vartheta_k(t)=\theta_k(t)-\Omega(\tau)t-\theta^0,~~~k=1,\ldots,N,
  \end{array}
\end{equation*}
setting $\theta^0=0$ for simplicity, we can rewrite equation~\eqref{eq:N-node_2nd-order_network}
with $x_1^{(j)}=\vartheta_j$, and  $x_2^{(j)}=\dot{\vartheta}_j$, as
\begin{equation}
  \label{eq:N_node_rotating_frame_vector_field}
  \begin{array}{l}
    \dot{x}_1^{(i)}=x_2^{(i)}\\
    \dot{x}_2^{(i)}=-\mu x_2^{(i)}-\mu\Omega(\tau)+\dfrac{K\mu}{N-1}\displaystyle\sum_{\stackrel[j\neq 1]{j=i}{}}^N\sin\left(x_{1\tau}^{(j)}-x_1^{(i)}-\Omega(\tau)\tau-\omega_M\tau \right),\\
 i=1,\ldots,N.
  \end{array}
\end{equation}
Then \eqref{eq:N_node_rotating_frame_vector_field} takes a form  similar to \eqref{e.f}.

In this section we  study bifurcations from $\S_N$-invariant relative equilibria. Hopf bifurcation from the
equilibrium $x=0$ in the comoving frame leads to  relative Hopf bifurcation of relative periodic orbits (RPOs).

\begin{definition}[Relative periodic orbit]
A relative  periodic orbit (RPO) of a $\Gamma$-equivariant dynamical system $\dot X = F(X)$ with phase space $\kX$ is a periodic orbit
in the space of group orbits $\kX/\Gamma$.
If $\Gamma$ is compact then a trajectory $X(t)$ lies on an RPO if and only if it is a periodic orbit in a comoving frame, see, e.g., ~\cite{Fiedler1996}.
\end{definition}
In our case $\Gamma = \sS^1 \times \S_N$ is compact and any RPO $\theta(t)$ satisfies
$\theta_j(t) = \theta_j(t-T) + C$ for all $j=1,\ldots, N$.
Then $\theta(t)$ is $T$-periodic in a frame moving
with velocity $C/T$.  

Linearising equation~\eqref{eq:N_node_rotating_frame_vector_field} around its equilibrium point $x=0$ we can
define  the matrix $L$  as in \eqref{eq:L(tau)_definition} 
\begin{equation*}
  L x=\left(
    \begin{array}{c}
      x_2^{(1)}\\
      -K\mu\cos(\widehat{\Omega}\tau)x_1^{(1)}-\mu x_2^{(1)}+\dfrac{K\mu}{N-1}\cos(\widehat{\Omega}\tau)e^{-\lambda\tau}\displaystyle\sum_{\stackrel[i\neq 1]{i=1}{}}^Nx_{1}^{(i)}\\

\vdots\\
      x_2^{(N)}\\
      -K\mu\cos(\widehat{\Omega}\tau)x_1^{(N)}-\mu x_2^{(N)}+\dfrac{K\mu}{N-1}\cos(\widehat{\Omega}\tau)e^{-\lambda\tau}\displaystyle\sum_{\stackrel[i\neq N]{i=1}{}}^Nx_{1}^{(i)}
    \end{array}
\right)
\end{equation*}
where 
\begin{equation}
\label{eq:Omega_omega_definition}
\widehat{\Omega}=\Omega(\tau)+\omega_M,
\end{equation}
or
\begin{equation}
  \label{eq:hat_omega_definition}
  \widehat{\Omega}=-K\sin(\widehat{\Omega}\tau)+\omega_M.
\end{equation}
From \eqref{eq:Omega(tau)_definition} we  see that
\begin{equation}
  \label{eq:Omega_omega_min_max}
  -K\leq \Omega(\tau) \leq K.
\end{equation}
Now, using the results obtained in
section~\ref{subsec:SN_symm_full_phase_model}, in particular
\eqref{eq:full_phase_linear_operator},~\eqref{eq:Nn_charact_functions_2}, and~
\eqref{eq:full_phase_P} with
\begin{equation*}
\begin{array}{rcl}
 q&=&K\mu\cos(\widehat{\Omega}\tau)\\
 r&=&-\dfrac{K\mu}{N-1}\cos(\widehat{\Omega}\tau)e^{-\lambda\tau},
\end{array}
\end{equation*}
 we obtain the characteristic functions
\begin{equation}
  \label{eq:Nn_charact_functions}
  \begin{array}{rcl}
P_{\Fix(\S_N)}(\lambda,\tau)&=&\text{det}(\triangle(\lambda,\tau)|_{\text{Fix}(\mathbf{S}_N)}))\\
&=&\lambda^2 + \mu\lambda + K\mu\cos(\widehat{\Omega}\tau) - K\mu\cos(\widehat{\Omega}\tau)e^{-\lambda\tau} \\
P_U(\lambda,\tau)&=& \text{det}(\triangle(\lambda,\tau)|_{V_j}), ~~j\neq 0\\
&=&\lambda^2 + \mu\lambda + K\mu\cos(\widehat{\Omega}\tau) + \dfrac{K\mu}{N-1}\cos(\widehat{\Omega}\tau)e^{-\lambda\tau}.
  \end{array}
\end{equation}
Clearly, $P_{\Fix(\S_N)}$ has a constant zero eigenvalue for all parameter values $\tau,K,\mu$ due to the translational symmetry. 
As before, roots in function $P_U$ correspond to eigenvalues of $L$ of  multiplicity $N-1$.

\begin{remark}({\em The parameter $k_v$ in the phase model and the rotating frame}) In~\cite{Correa2013} a modification of the model~\eqref{eq:N-node_2nd-order_network} was presented by introducing the parameter $k_v$ in order to avoid a zero eigenvalue in the characteristic equation $P_{\Fix(\S_N)}(\lambda,\tau)=0$ in~\eqref{eq:Nn_charact_functions}; the phase model
from~\cite{Correa2013} is
\begin{equation}
  \ddot{\theta}_i(t)+(\mu+k_v)\dot{\theta}_i(t)+\mu k_v\theta_i(t)-\frac{K\mu}{N-1}\sum_{\stackrel[j\neq i]{j=1}{}}^N\sin(\theta_j(t-\tau)-\theta_i(t)-\omega_M\tau)=0,
\end{equation}
for $ i = 1,\ldots, N.$
  The $\R$-symmetry  disappears when $k_v\neq0$. 
\end{remark}


\subsection{Symmetry-preserving bifurcations}
\label{subsec:Bif_analysis_phase_model}
The rotation
frequency $\Omega(\tau)$ of the  $\S_N$-invariant relative equilibria from 
\eqref{e.RE} is determined by \eqref{eq:Omega(tau)_definition}; for a given $\tau$ there
exists a whole family of solutions $\Omega(\tau)$ satisfying this
 equation, and,  as $\tau$ increases, more solutions appear. We now fix study symmetry preserving bifurcation from $\S_N$-invariant relative equilibria. 
\marginnote{
\tiny
\begin{align*}
R(i\omega,\tau)&=-\omega^2+K\mu\cos(\widehat{\Omega}\tau)+i\mu\omega
\end{align*}
\normalsize
}
From ~\eqref{eq:Nn_charact_functions} when $\tau=0$ we have
\begin{equation}
  P_{\Fix(\S_N)}(\lambda,0)=\lambda^2+\mu\lambda,
\end{equation}
whose roots are
\begin{equation}
  \lambda_+=0, ~~~\lambda_-=-\mu.\nonumber
\end{equation}
\marginnote{
\tiny
\begin{align*}
  \sin(\omega\tau)&=\dfrac{R_IS_R-S_IR_R}{|S|^2} = \frac{R_I}{S_R} = 
\\                &=-\dfrac{\omega}{K\cos(\widehat{\Omega}\tau)}
\\\cos(\omega\tau)&=-\dfrac{S_IR_I+S_RR_R}{|S|^2} =  -\frac{R_R}{S_R} =
\\                &=\dfrac{-\omega^2+K\mu\cos(\widehat{\Omega}\tau)}{K\mu\cos(\widehat{\Omega}\tau)}
\end{align*}
\normalsize
}
In order to find critical delays leading to bifurcations in
$\Fix(\S_N)$ we will follow section~\ref{subsec:Sn_map}. Using~\eqref{eq:Sn_map_P} we have that $P_{\Fix(\S_N)}(\lambda,\tau)=R(\lambda,\tau)+S(\tau)e^{-\lambda\tau}$ where
\begin{equation}
  \label{eq:phase_model_R_S_Fix}
  \begin{array}{rcl}
    R(\lambda,\tau)&=&\lambda^2+\mu\lambda+K\mu\cos(\widehat{\Omega}\tau)\\\\
    S(\tau)&=&-K\mu\cos(\widehat{\Omega}\tau),
  \end{array}
\end{equation}
and from \eqref{eq:Sn_map_sin_cos} we obtain
\marginnote{
\tiny
\begin{align*}
F(\omega,\tau)&=R(i\omega,\tau)R(-i\omega,\tau)
\\            &~-S(i\omega,\tau)S(-i\omega,\tau)=0
\\            &=|R(i\omega,\tau)|^2-|S(\tau)|^2=0
\\            &=\left(-\omega^2+K\mu\cos(\widehat{\Omega}\tau)\right)^2
\\            &~+(\mu\omega)^2-\left(K\mu\cos(\widehat{\Omega}\tau)\right)^2=0 
\\            &=\omega^4 - 2K\mu\cos()\omega^2 + K^2\mu^2\cos^2() 
\\            &~+\mu^2\omega^2- K^2\mu^2\cos^2()
\end{align*}
\normalsize
}
\begin{equation}
  \label{eq:2n_P1_sin_cos_a}
  \begin{array}{lll}
    \sin(\omega\tau)&=&-\dfrac{\omega}{K\cos(\widehat{\Omega}\tau)}\\\\
    \cos(\omega\tau)&=&\dfrac{-\omega^2+K\mu\cos(\widehat{\Omega}\tau)}{K\mu\cos(\widehat{\Omega}\tau)}
  \end{array},
\end{equation}
for $K\mu\cos(\widehat{\Omega}\tau)\neq0$.
\marginnote{
\tiny
\begin{align*}
\cos(\widehat{\Omega}\tau)&= \pm \sqrt{ 1 - \sin^2(\widehat{\Omega}\tau)}
\\ & =  \pm \sqrt{ 1 -  \frac{(\omega_M - \hat\Omega)^2}{K^2} }
\end{align*}
\normalsize
}
The polynomial $F(\omega,\tau)$ from~\eqref{eq:Sn_map_F} becomes
\begin{equation}
  \begin{array}{rcl}
  F(\omega,\tau)&=&\omega^2(\omega^2-2K\mu\cos(\widehat{\Omega}\tau)+\mu^2)
  \end{array}
\end{equation}
with roots
\begin{equation}
  \label{eq:2n_omega}
  \omega^2={2K\mu\cos(\widehat{\Omega}\tau)-\mu^2},~~~\text{or}~~~\omega=0.
\end{equation}
From \eqref{eq:hat_omega_definition} we obtain
\begin{equation}
  \cos(\widehat{\Omega}\tau)=\pm\dfrac{1}{K}\sqrt{K^2-(\omega_M-\widehat{\Omega})^2},
\end{equation}
\marginnote{
\tiny
\begin{align*}
\widehat{\Omega}'_{\tau}=&~\dfrac{d\widehat{\Omega}}{d\tau}
\\=&~-K\cos(\widehat{\Omega}\tau)(\widehat{\Omega}+\tau\widehat{\Omega}'_{\tau})
\\\widehat{\Omega}'_{\tau}=&~-\dfrac{\widehat{\Omega}K\cos(\widehat{\Omega}\tau)}{1+\tau K\cos(\widehat{\Omega}\tau)}.  
\end{align*}
\normalsize
}
thus 
\begin{equation}
  \label{eq:2n_omega_a}
  \omega=\pm\left({\pm 2\mu\sqrt{K^2-(\omega_M-\widehat{\Omega})^2} -\mu^2}\right)^{1/2},~~~\text{or}~~~\omega=0.
\end{equation}
\marginnote{ \tiny
\begin{align*}
 A&=\Re\left((\cos(\omega\tau)-\i\sin(\omega\tau))\right.\\
& ~~~\left. \left(-\i\omega K\mu\cos(\widehat\Omega\tau)-K\mu\sin(\widehat\Omega\tau)\right.\right.\\
&~~~\left.\left. \left(\widehat\Omega+\tau\widehat\Omega'_\tau \right)\right)\right.
\\ &~~~\left.+K\mu\sin(\widehat\Omega\tau)\left(\widehat\Omega+\tau\widehat\Omega'_\tau \right)\right)
\\ &=-\cos(\omega\tau)K\mu\sin(\widehat\Omega\tau)\\
&~~~ \left(\widehat\Omega+\tau\widehat\Omega'_\tau \right)-\\
&~~~ \left. \sin(\omega\tau)\omega K\mu\cos(\widehat\Omega\tau) \right.\\
& ~~~ \left.+K\mu\sin(\widehat\Omega\tau)\left(\widehat\Omega+\tau\widehat\Omega'_\tau \right)\right.\\
\\ &=(\omega^2-K\mu\cos(\widehat\Omega\tau))\frac{\sin(\widehat\Omega\tau)}{\cos(\widehat\Omega\tau)}\\
& ~~~\left(\frac{\widehat\Omega}{1+\tau K\cos(\widehat\Omega\tau)}\right)\\
& ~~~+\omega^2\mu+K\mu\sin(\widehat\Omega\tau)\frac{\widehat\Omega}{1+\tau K\cos(\widehat\Omega\tau)}
\\ &=\frac{\omega^2\widehat\Omega\sin(\widehat\Omega\tau)}{\cos(\widehat\Omega\tau)(1+\tau K\cos(\widehat\Omega\tau)) }+\omega^2\mu. \end{align*}
}
Solutions $\omega\in\R^+$ exist provided
\begin{equation}
  \label{eq:2n_con_1a}
   \cos(\widehat{\Omega}\tau)\geq 0~~\mbox{and}~~ 2\sqrt{K^2-(\omega_M-\widehat{\Omega})^2}\geq\mu.
   \end{equation}
Given $\tau\in\R^+$, we can compute $\widehat{\Omega}$ using~\eqref{eq:Omega_omega_definition} and~\eqref{eq:Omega(tau)_definition}. For $\omega\in\R^+$  satisfying \eqref{eq:2n_omega_a} we compute the $S_n$ map, see section~\ref{subsec:Sn_map}, whose zeros are the critical bifurcation time delays for $\Fix(\S_N)$.
Using \eqref{eq:Sn_delta_abcd} we obtain $\delta(\omega(\tau^*))$  from \eqref{eq:2n_delta} to find the direction in which roots, 
if any, cross the imaginary axis. From \eqref{eq:hat_omega_definition} we compute
\begin{equation}
  \label{eq:dOmega_omega_dtau}
  \begin{array}{c} 
  \widehat{\Omega}'_{\tau}=\dfrac{d\widehat{\Omega}}{d\tau}=-\dfrac{\widehat{\Omega}K\cos(\widehat{\Omega}\tau)}{1+\tau K\cos(\widehat{\Omega}\tau)}.
  \end{array}
\end{equation}
\marginnote{\tiny
\begin{align*}
B&=\Im\left((\cos(\omega\tau)-\i\sin(\omega\tau))\right.
\\ &~~~\left.\left(-\i\omega K\mu\cos(\widehat\Omega\tau) \right.\right.
\\
 &~~~\left.\left.-K\mu\sin(\widehat\Omega\tau)\left(\widehat\Omega+\tau\widehat\Omega'_\tau \right)\right)\right.
\\ &~~~\left.+K\mu\sin(\widehat\Omega\tau)\left(\widehat\Omega+\tau\widehat\Omega'_\tau \right)\right)
\\ &=-\cos(\omega\tau)\omega K \mu \cos(\widehat\Omega\tau)+ \sin(\omega\tau)K \mu
\\ &~~~* \sin(\widehat\Omega\tau)(\widehat\Omega+\tau\widehat\Omega'_\tau)
\\ &=\left(\omega^2-K\mu\cos(\widehat\Omega\tau) \right)\omega
\\ &~~~
-\frac{\omega\mu}{\cos(\widehat\Omega\tau)}  \frac{\widehat\Omega\sin(\widehat\Omega\tau)}{(1+\tau K\cos(\widehat\Omega\tau))}
\\
C&= \Re( R_\lambda + e^{-i\omega \tau}(S_\lambda-\tau S))\\
&= \Re( 2i\omega + \mu -(\cos(\omega\tau)-i\sin(\omega\tau))*
\\ & *K\mu\tau \cos(\hat\Omega \tau)\\
& = \mu - \cos(\omega\tau)K\mu\tau \cos(\hat\Omega \tau)\end{align*}
\normalsize
}
and, from~\eqref{eq:phase_model_R_S_Fix},~\eqref{eq:2n_P1_sin_cos_a}  
and \eqref{eq:Sn_map_abcd}
\begin{align}\label{e.ABCDPhaseModel}
\begin{array}{l}
    A= \dfrac{\omega^2\widehat\Omega\sin(\widehat\Omega\tau)}{\cos(\widehat\Omega\tau)(1+\tau K\cos(\widehat\Omega\tau)) }+\omega^2\mu\\
   B= \omega\left(\omega^2-K\mu\cos(\widehat\Omega\tau)-\dfrac{\mu\widehat\Omega\sin(\widehat\Omega\tau)}{\cos(\widehat\Omega\tau)(1+\tau K\cos(\widehat\Omega\tau))} \right)\\
    C= \mu+\tau\left(K\mu\cos(\widehat{\Omega}\tau)-\omega^2\right)\\
    D= \omega(2+\tau\mu).
\end{array}
 \end{align}


\subsubsection{Symmetry preserving bifurcation  from equilibria}\label{sss:symPreservFromEquil}
\marginnote{\tiny
\begin{align*}
D&= \Im() = 2\omega -\sin(\omega\tau)K\mu \tau\cos(\widehat{\Omega}\tau)
\end{align*}
\normalsize
}
From \eqref{eq:Omega(tau)_definition} we see that $\Omega(\tau)=0$ is a rotating co-frame solution and hence
relative equilibrium becomes an equilibrium 
when $\omega_M\tau=n\pi$, with $n\in\N_0$. Now, here we have two possible cases:
\marginnote{ 
\tiny  
\begin{align*}
  R(\lambda;\tau)&=\lambda^2+\mu\lambda+K\mu
\\R(\i\omega)&=-\omega^2+K\mu+\i\mu\omega
\\S(\lambda;\tau)&=-K\mu. 
\end{align*}
  \begin{align*}
    \sin(\omega\tau)&=\frac{R_IS_R-S_IR_R}{|S|^2}
\\                  &=-\frac{\omega}{K}
\\  \cos(\omega\tau)&=-\frac{S_IR_I+S_RR_R}{|S|^2}
\\                  &=\frac{K\mu-\omega^2}{K\mu}.
  \end{align*}
\begin{align*}
\hat\Omega + \tau \hat\Omega_\tau=\frac{\hat\Omega}{1+\tau K \cos(\hat\Omega\tau)}
\end{align*}
}
\begin{itemize}
\item When $\omega_M\tau=2n\pi$. We have $\widehat{\Omega}=\omega_M$ and $\cos(\widehat{\Omega}\tau)=1$, thus from \eqref{eq:2n_P1_sin_cos_a} we obtain
\begin{equation}
\begin{array}{rcl}
 \sin(\omega\tau)&=&-\dfrac{\omega}{K}\\
 \cos(\omega\tau)&=&\dfrac{K\mu-\omega^2}{K\mu},
\end{array}
\label{eq:sin_cos_tau_2npi}
\end{equation}
and from condition~\eqref{eq:2n_omega} we obtain 
\begin{equation}
  \label{eq:omega_tau_2npi}
  \begin{array}{rcl}
  \omega&=&\pm\sqrt{2K\mu-\mu^2},
  \end{array}
\end{equation}
provided $2K\geq\mu$.
Then, from~\eqref{eq:sin_cos_tau_2npi} and~\eqref{eq:omega_tau_2npi}
we obtain a second condition
\marginnote{
\tiny
\begin{align*}
   \sin(\omega\tau)=&~\dfrac{R_IS_R-S_IR_R}{|S|^2}
\\                 =&~\dfrac{R_I}{S_R}
\\                 =&~\dfrac{\omega\mu}{-K\mu\cos(\widehat{\Omega}\tau)}
\\                 =&~-\dfrac{\omega}{K}
\\ \cos(\omega\tau)=&~-\dfrac{S_IR_I+S_RR_R}{|S|^2}
\\                 =&~-\dfrac{R_R}{S_R}
\\                 =&~-\dfrac{-\omega^2+K\mu\cos(\widehat{\Omega}\tau)}{-K\mu\cos(\widehat{\Omega}\tau)}
\\                 =&~\dfrac{K\mu-\omega^2}{K\mu}.
\\\text{\eqref{eq:omega_tau_2npi}}~\omega=&~\pm\sqrt{2K\mu\cos(\widehat{\Omega}\tau)-\mu^2 }.
\end{align*}
\normalsize
}
\begin{equation}
\label{eq:omega_tau_2npi_2nd_cond}
\sqrt{2K\mu-\mu^2} 
=~\dfrac{\omega_M}{2n\pi}\left(\arg\left( \mu -  K,-\sqrt{2K\mu-\mu^2 }\right)+2m\pi \right),
\end{equation}
with $n\in\N,~m\in\Z$.

The solution of \eqref{eq:omega_tau_2npi_2nd_cond} is  a curve 
$K(\mu;\omega_M,n,m)$ for which a   critical delay at $\omega_M\tau=2n\pi$ exists with imaginary eigenvalue $\lambda=\pm \i\omega(K,\mu)$ with $\omega$ from~\eqref{eq:omega_tau_2npi}, and $\widehat{\Omega}=\omega_M$. 
 In figure~\ref{fig:curves_K_pm_phase_model} the curves $K(\mu;\omega_M,n,m)$   are shown
 for $m=\{1,\ldots,4\}$, $\omega_M=1$ and $n=1$.
 Now, calculating $\delta(\omega,\tau)$ for the case $\tau=2n\pi/\omega_M$ 
using~\eqref{eq:Sn_delta_abcd} and~\eqref{e.ABCDPhaseModel} we compute 
the denominator $AC+BD$ which determines the sign of $\delta(\omega,\tau)$ as
 in \eqref{eq:full_phase_delta_abcd_general} with $\cos(2\phi^*)=0$.
  As before, \eqref{eq:full_phase_AC+BD_b} holds true
 with $\omega=\omega_+$.
 From these curves  
families of periodic orbits bifurcate.
\begin{figure}[!htb]
\centering
\includegraphics[scale=0.4]{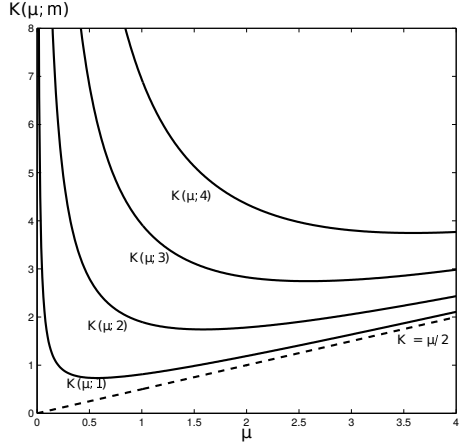}
\caption{Curves $K(\mu;\omega_M,n,m)$
  from~\eqref{eq:omega_tau_2npi_2nd_cond} with $\omega_M=1$, $n=1$  for different values of $m$.}
\label{fig:curves_K_pm_phase_model}
\end{figure}
\item When $\omega_M\tau=(2n+1)\pi$
 there are no non-zero real roots
of \eqref{eq:2n_omega}.
\end{itemize}
\marginnote{\tiny
  \begin{align*}
 &   \omega\tau=\Arctan\left(\frac{-\omega\mu}{K\mu-\omega^2}\right)\\
&+2m\pi
\\
& \pm\sqrt{2K\mu-\mu^2}\left(\frac{2n\pi}{\omega_M}\right) 
\\ & = \Arctan\left(\frac{\mp\mu\sqrt{2K\mu-\mu^2}}{K\mu-(2K\mu-\mu^2
    )}\right)+2m\pi            
  \end{align*}
}

\subsubsection{Symmetry preserving bifurcation from relative equilibria}
\label{subsubsec:Fix_SN_Omega_neq_0}
\begin{figure}[!htb]
  \centering
  \includegraphics[scale=.5]{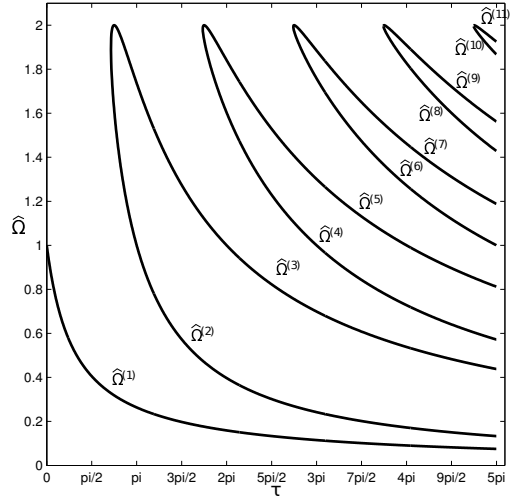}
  \caption{Curves of relative equilibria determined by their rotation frequency $\widehat{\Omega}$, for $\mu=1$, $K=1$, and $\omega_M=1$.}
  \label{fig:Omega_w_curves}
\end{figure}

\begin{figure}[!htb]
\centering
\includegraphics[scale=.45]{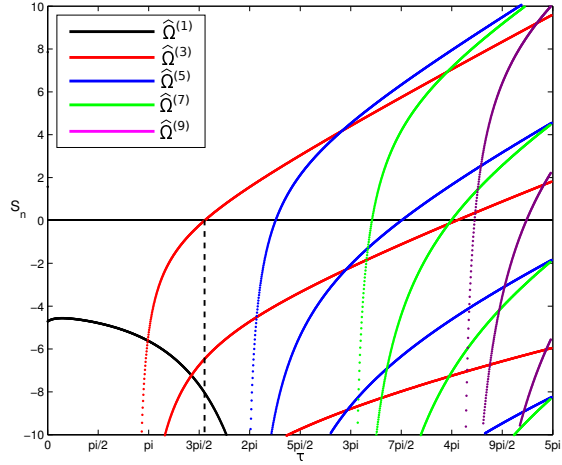}
\caption{The $S_n$ map for several curves of relative equilibria  determined by their rotation frequency $\widehat{\Omega}^{(j)}(\tau)$.   The zeroes of the $S_n$ map
  determine imaginary roots of $P_{\Fix(\S_N)}(\lambda,\tau)$ from \eqref{eq:Nn_charact_functions}.
 }
 \label{fig:Sn_map_Pfix}
\end{figure}

For this case, when $\Omega(\tau) \neq 0$, we carried out numerical computations to find 
 time-delays $\tau$ leading to bifurcations in $\Fix(\S_N)$; the procedure is as follows:
\begin{itemize}
\item For a given $\tau\in\R^+$ and parameters $K,\omega_M,\mu\in\R^+$ we calculate all   $m$ real solutions $\widehat{\Omega}^{(j)}=\widehat{\Omega}^{(j)}(\tau;K,\omega_M,\mu)$, $j=1,\ldots,m$ of~\eqref{eq:hat_omega_definition}, to determine all $\S_N$-invariant relative equilibria,  noting \eqref{eq:Omega_omega_min_max} and that the number of solutions is finite and increases with $\tau$.
%
\marginnote{
\tiny
\begin{align*}
\dfrac{\sin(\omega\tau)}{\cos(\omega\tau)}&=  -\tfrac{\omega}{K} \tfrac{K\mu}{K\mu - \omega^2} =
-\tfrac{\omega\mu}{K\mu - \omega^2} \\
&= -\tfrac{\omega\mu}{K\mu - (2K\mu-\mu^2)}   = \tfrac{\omega}{\mu -  K} \\
&=~\dfrac{\pm\sqrt{2K\mu-\mu^2}}{\mu -  K}.
\end{align*}
and
\[
\omega \tau = \Arctan\left(\dfrac{\pm\sqrt{2K\mu-\mu^2}}{\mu -  K}\right)
\]
\normalsize
}
\item For each solution $\widehat{\Omega}^{(j)}$ we compute $\omega=\omega_{\pm}$ from~\eqref{eq:2n_omega_a}, provided condition~\eqref{eq:2n_con_1a} holds. 
\item For each 
$\omega=\omega\left(\widehat{\Omega}^{(j)};K,\omega_M,\mu\right)\in\R^+$ we compute the $S_n$ map using \eqref{eq:2n_P1_sin_cos_a}, see section~\ref{subsec:Sn_map}.
 The $S_n$ map gives us the time-delay $\tau_n(\tau)$ from \eqref{eq:2n_tau_map}, which depends on the values calculated previously. If $\tau_n(\tau)$  matches   the given $\tau$, then we have find a bifurcation time-delay, cf.~\eqref{eq:2n_Sn_map}.
\item Finally in order to determine the direction in which these roots cross the imaginary axis we have to compute the sign of $\delta(\omega(\tau^*))$ using \eqref{eq:Sn_delta_abcd}.
\end{itemize}
In figure~\ref{fig:Omega_w_curves} the eleven possible curves $\widehat{\Omega}^{(j)}$ within the interval $\tau\in[0,5\pi]$ are shown with parameters $\mu=1$, $K=1$ and $\omega_M=1$.
In figure~\ref{fig:Sn_map_Pfix} the $S_n$ maps for those
curves of relative equilibria within the interval $\tau\in[0,5\pi]$  
 are shown. At $S_n=0$ along the curves $\widehat{\Omega}^{(j)}$  relative Hopf bifurcation of synchronized RPOs occurs.
The sign for $\delta(\omega(\tau^*))$ is positive in all cases, i.e.,  the roots cross the imaginary axis from the left to the right.


\subsection{Symmetry breaking bifurcations}
\label{subsec:bif_Xj_phase_model}
When $\tau=0$ the characteristic function $P_U$ from \eqref{eq:Nn_charact_functions}  becomes
\begin{equation*}
  P_U(\lambda,0)=\lambda^2+\mu\lambda+K\mu\left(\frac{N}{N-1}\right),
\end{equation*}
whose roots are
\begin{equation*}
  \lambda_{\pm}=-\frac{\mu}{2}\pm\frac{1}{2}\left(\mu^2-4K\mu\frac{N}{N-1} \right)^{1/2},
\end{equation*}
and since $\mu,K\in\R^+$ and $N\geq2$, nonzero roots of $P_U$ when $\tau=0$ are always stable.
%
\marginnote{
\tiny
\begin{align*}
    \sin(\omega\tau)=&~\dfrac{R_IS_R-S_IR_R }{|S|^2}
\\                  =&~\dfrac{R_I}{S_R}
\\                  =&~\dfrac{\omega(N-1)}{K\cos(\widehat{\Omega}\tau)}
\\  \cos(\omega\tau)=&~-\dfrac{S_IR_I+S_RR_R}{|S|^2}
\\                  =&~-\dfrac{R_R}{S_R}
\\                  =&~\dfrac{\left(\omega^2-K\mu\cos(\widehat{\Omega}\tau)\right)(N-1)}{K\mu\cos(\widehat{\Omega}\tau)}
\\&
\\&
\\F(\omega,\tau)=&~R_R^2+R_I^2-S_R^2 
\\              =&~\left(-\omega^2+K\mu\cos(\widehat{\Omega}\tau)\right)^2
\\               &~+(\mu\omega)^2-\left(\dfrac{K\mu}{N-1}\cos(\widehat{\Omega}\tau) \right)^2
\\              =&~\omega^4+(\mu^2-2K\mu\cos(\widehat{\Omega}\tau))\omega^2
\\               &~+(1-\dfrac{1}{(N-1)^2})K^2\mu^2\cos^2(\widehat{\Omega}\tau)
\\              =&~\omega^4+(\mu^2-2K\mu\cos(\widehat{\Omega}\tau))\omega^2
\\               &~+\dfrac{N(N-2)}{(N-1)^2}K^2\mu^2\cos^2(\widehat{\Omega}\tau)
\end{align*}
\normalsize
}
In order to calculate the $S_n$ map, see section~\ref{subsec:Sn_map}, in particular \eqref{eq:Sn_map_P},
 we note that by \eqref{eq:Nn_charact_functions},
\begin{equation}
  \label{eq:phase_model_Xj_R_S}
  \begin{array}{rcl}
    R(\lambda,\tau)&=&\lambda^2+\mu\lambda+K\mu\cos(\widehat{\Omega}\tau)\\
    S(\tau)&=&\dfrac{K\mu}{N-1}\cos(\widehat{\Omega}\tau),
  \end{array}
\end{equation}
so from \eqref{eq:Sn_map_sin_cos} we obtain
\begin{equation}
  \label{eq:Nn_sin_cos}
  \begin{array}{rcl}
    \sin(\omega\tau)&=&\dfrac{\omega(N-1)}{K\cos(\widehat{\Omega}\tau)}\\
    \cos(\omega\tau)&=&\dfrac{\left(\omega^2-K\mu\cos(\widehat{\Omega}\tau)\right)(N-1)}{K\mu\cos(\widehat{\Omega}\tau)},
  \end{array}
\end{equation}
and   the polynomial equation from \eqref{eq:Sn_map_F} is
\begin{align}
F(\omega,\tau)=&~\omega^4+(\mu^2-2K\mu\cos(\widehat{\Omega}\tau))\omega^2+\dfrac{N(N-2)}{(N-1)^2}K^2\mu^2\cos^2(\widehat{\Omega}\tau)=0,
\end{align}
here again, $\widehat{\Omega}=\Omega(\tau)+\omega_M$ and $\Omega(\tau)$ satisfies  \eqref{eq:Omega(tau)_definition}.
Hence,
\begin{equation}
  \label{eq:Nn_omega}
\begin{array}{lll}
  \omega^2_{\pm}&=&-\dfrac{1}{2}\left(\mu^2-2K\mu\cos(\widehat{\Omega}\tau)\right)\\
&&\pm\dfrac{1}{2}\left[\left(\mu^2-2K\mu\cos(\widehat{\Omega}\tau)  \right)^2-4\dfrac{N(N-2)}{(N-1)^2}K^2\mu^2\cos^2(\widehat{\Omega}\tau)\right]^{1/2}.
\end{array}
\end{equation}
\marginnote{
\tiny
\begin{align*}
\left(\mu^2-2K\mu\cos(\widehat{\Omega}\tau)\right)^2\geq&
\\4\dfrac{N(N-2)}{(N-1)^2}K^2\mu^2\cos^2(\widehat{\Omega}\tau)&
\end{align*}
\normalsize
}
In~\eqref{eq:Nn_omega} the discriminant is always smaller than the square of the first term, therefore in order to $\omega\in\R$ this first term has to be non-negative, i.e.,
\begin{equation}
  \mu\leq 2K\cos(\widehat{\Omega}\tau).
\end{equation}
Hence, $\cos(\widehat{\Omega}\tau)\geq 0$.
For $\omega \in \R$ the discriminant also has to be greater or equal to zero, i.e.,
\marginnote{\tiny
\begin{align*}
  A&=\Re\Bigg[(\cos(\omega\tau)-\i\sin(\omega\tau))\big[\i\omega\frac{K\mu}{N-1}\cos(\widehat\Omega\tau)\\
&~~~+\frac{K\mu}{N-1}\sin(\widehat\Omega\tau)(\widehat\Omega+\tau\widehat\Omega'_\tau) \big]
\\ &~~~+K\mu\sin(\widehat\Omega\tau)(\widehat\Omega+\tau\widehat\Omega'_\tau) \Bigg]
\\ &=\cos(\omega\tau)\frac{K\mu}{N-1}\sin(\widehat\Omega\tau)(\widehat\Omega+\tau\widehat\Omega'_\tau)\\
&~~~+\sin(\omega\tau)\omega\frac{K\mu}{N-1}\cos(\widehat\Omega\tau)\\
&~~~+K\mu\sin(\widehat\Omega\tau)(\widehat\Omega+\tau\widehat\Omega'_\tau)
\\ &=(\omega^2-K\mu\cos(\widehat\Omega\tau))
\\
&~~~\frac{\sin(\widehat\Omega\tau)\widehat\Omega}{\cos(\widehat\Omega\tau)(1+\tau K\cos(\widehat\Omega\tau))}\\
&~~~+\omega^2\mu+\frac{K\mu\sin(\widehat\Omega\tau)\widehat\Omega}{1+\tau K\cos(\widehat\Omega\tau)}
\\ &=\frac{\omega^2\sin(\widehat\Omega\tau)\widehat\Omega}{\cos(\widehat\Omega\tau)(1+\tau K\cos(\widehat\Omega\tau))}+\omega^2\mu
\\ &=\omega^2\left(\frac{\sin(\widehat\Omega\tau)\widehat\Omega}{\cos(\widehat\Omega\tau)(1+\tau K\cos(\widehat\Omega\tau))}+\mu\right)
\end{align*}
\normalsize}
\begin{equation}
\label{eq:phase_model_omega_condition}
\ 2K\cos(\widehat{\Omega}\tau) - \mu \geq
\dfrac{\sqrt{N(N-2)}}{N-1}2K\cos(\widehat{\Omega}\tau).
\end{equation}
from which we get:
\begin{equation}
 \mu\leq 2K\cos(\widehat{\Omega}\tau)\left(1-\dfrac{\sqrt{N(N-2)}}{N-1}\right).
  \label{eq:phase_model__condition_b}
\end{equation}
Using~\eqref{eq:2n_delta} and~\eqref{eq:Sn_delta_abcd} we obtain $\delta(\omega(\tau^*))$ to find the direction in which roots, 
if any, cross the imaginary axis. From \eqref{eq:Sn_map_abcd} and
noting~\eqref{eq:dOmega_omega_dtau},~\eqref{eq:phase_model_Xj_R_S},HERE!!!!!
and~\eqref{eq:Nn_sin_cos} we see that $A$, $B$, $C$ and $D$ are as in \eqref{e.ABCDPhaseModel}.


\subsubsection{Symmetry-breaking bifurcations of relative equilibria}
Setting $\lambda=0$  \eqref{eq:Nn_charact_functions} gives   time delays $\tau=\tau_*$  
\begin{align}
\label{eq:Omega_tau=pi/2_a}
\widehat{\Omega}\tau^*=\dfrac{\pi}{2}+n\pi,~~~n\in\Z,  
\end{align}
where $P_U$ has a zero root.
\marginnote{\tiny
\begin{align*}
\\B&=\Im\Bigg[(\cos(\omega\tau)-\i\sin(\omega\tau))\big[\i\omega\frac{K\mu}{N-1}\cos(\widehat\Omega\tau)\\
&~~~+\frac{K\mu}{N-1}\sin(\widehat\Omega\tau)(\widehat\Omega+\tau\widehat\Omega'_\tau) \big]
\\ &~~~+K\mu\sin(\widehat\Omega\tau)(\widehat\Omega+\tau\widehat\Omega'_\tau) \Bigg]
\\ &=\cos(\omega\tau)\frac{\omega K\mu}{N-1}\cos(\widehat\Omega\tau)
\\&~~~-\sin(\omega\tau)\frac{K\mu}{N-1}\sin(\widehat\Omega\tau)(\widehat\Omega+\tau\widehat\Omega'_\tau)
\\ &=(\omega^2-K\mu\cos(\widehat\Omega\tau))\omega-
\\ &~~~\omega\mu\frac{\sin(\widehat\Omega\tau)}{\cos(\widehat\Omega\tau)}\left(\frac{\widehat\Omega}{1+\tau K\cos(\widehat\Omega\tau)}\right)
\\ &=\omega\left(\omega^2-K\mu\cos(\widehat\Omega\tau)\right.
\\ &~~~-\left.\frac{\mu\widehat\Omega\sin(\widehat\Omega\tau)}{\cos(\widehat\Omega\tau)(1+\tau K\cos(\widehat\Omega\tau))} \right)
\end{align*}
\normalsize
}
From \eqref{eq:hat_omega_definition} we see that
\begin{align}
\widehat\Omega=(-1)^{n+1} K+\omega_M.
\end{align}
Substituting this into~\eqref{eq:Omega_tau=pi/2_a} we obtain 
\begin{align}
  \label{eq:tau_critical_b}
  \tau^*= 
    \left(\dfrac{\pi}{2}+n\pi\right)\dfrac{1}{\omega_M+ (-1)^{n+1} K },~~~n\in\Z,
\end{align}
where $n\in \Z$ such that $\tau^*\geq 0$.
Now we can calculate $\delta(\omega,\tau)|_{\omega=0}$
using~\eqref{eq:Sn_delta_abcd}, \eqref{eq:Nn_charact_functions} and \eqref{eq:dOmega_omega_dtau},
\begin{align}
  \delta(0,\tau^*)= 
      \dfrac{(-1)^{n} KN}{N-1}(\omega_M+ (-1)^{n+1} K),
\end{align}
where $n \in \Z$ is as in \eqref{eq:tau_critical_b}. At these critical time delays relative equilibria bifurcate which are not $\S_N$-invariant,
i.e., they satisfy $\theta_k(t) = \Omega(\tau) t + \theta_k^0$ with $\theta_k^0 \neq \theta_j^0$ for some
$j\neq k$, for details see \cite{Golubitsky1988},  cf. also Section \ref{subsec:spatio-temp_symm}.
\marginnote{ \tiny
\begin{align*}
  C&=\Re\bigg[\mu+2\i\omega+(\cos(\omega\tau)
\\ &~~-\i\sin(\omega\tau))\left(-\tau\frac{K\mu}{N-1}\cos(\widehat\Omega\tau)\right) \bigg]
\\ &=\mu-\cos(\omega\tau)\tau\frac{K\mu}{N-1}\cos(\widehat\Omega\tau)
\\ &=\mu-\tau(\omega^2-K\mu\cos(\widehat\Omega\tau))
\\D&=\Im\bigg[\mu+2\i\omega+(\cos(\omega\tau)
\\ &~~-\i\sin(\omega\tau))\left(-\tau\frac{K\mu}{N-1}\cos(\widehat\Omega\tau)\right) \bigg]
\\ &=2\omega+\sin(\omega\tau)\frac{\tau K\mu}{N-1}\cos(\widehat\Omega\tau)
\\ &=\omega(2+\tau\mu). 
\end{align*}
}

\subsubsection{Symmetry-breaking bifurcations from equilibria}
We shall analyze Hopf bifurcation from equilibria where $\Omega(\tau)=0$ as we did in  Section  \ref{sss:symPreservFromEquil} for $\Fix(\S_N)$. We know that for this case $\omega_M\tau=n\pi$ with $n\in\N_0$ and $\widehat{\Omega}=\omega_M$; we have two cases:
\begin{itemize}
\item 
When $\omega_M\tau=2n\pi$ then $\cos(\omega_M\tau)=1$  and Hopf bifurcations with frequencies $\omega_\pm$ given by \eqref{eq:Nn_omega} are possible
if \eqref{eq:phase_model__condition_b} holds with $\cos(\widehat{\Omega}\tau)=1$.
Now, calculating $\delta(\omega,\tau)$ for the case $\tau=2n\pi/\omega_M$ 
using~\eqref{eq:Sn_delta_abcd} and~\eqref{e.ABCDPhaseModel} we compute 
the denominator $AC+BD$ which determines the sign of $\delta(\omega_\pm,\tau)$  
as in \eqref{eq:full_phase_delta_abcd_general} with $\cos(2\phi^*)=0$.
Hence, \eqref{eq:full_phase_AC+BD_b} holds true again.
In this case equivariant Hopf bifurcation takes place and families of non-synchronous periodic orbits bifurcate, see Section \ref{subsec:spatio-temp_symm}.
\item When $\omega_M\tau=(2n+1)\pi$. We have $\cos(\omega_M\tau)=-1$
which violates \eqref{eq:phase_model__condition_b}, therefore symmetry breaking bifurcations are not possible in  this case.
\end{itemize}
\marginnote{\tiny 
  \begin{align*}
   & \dfrac{\partial P_j}{\partial\tau}=-K\mu\sin(\widehat\Omega\tau)\left(\widehat\Omega+\tau\widehat\Omega'_{\tau} \right)\\
&~~-\dfrac{K\mu}{N-1}\cos(\widehat\Omega\tau)\lambda e^{-\lambda\tau}
\\&~~~-e^{-\lambda\tau}\dfrac{K\mu}{N-1}\sin(\widehat\Omega\tau)\left(\widehat\Omega+\tau\widehat\Omega'_{\tau} \right)
\\   & \dfrac{\partial P_j}{\partial\tau}(0,\tau)|_{\widehat\Omega\tau=\frac{\pi}{2}+n\pi}=(-1)^{n+1} \frac{K\mu\widehat\Omega N}{N-1},
  \end{align*}
using that $\widehat\Omega'_{\tau}|_{\widehat\Omega\tau=\frac{\pi}{2}+n\pi}=0$.
\begin{align*}
  \dfrac{\partial P_j}{\partial\lambda}(0,\tau)&=2\lambda+\mu-\dfrac{K\mu}{N-1}\cos(\widehat\Omega\tau)\tau e^{-\lambda\tau}
\\\dfrac{\partial P_j}{\partial\lambda}(0,\tau)&=\mu.
\end{align*}
}
\subsubsection{Curves of symmetry-breaking bifurcations}
\begin{figure}[!htb]
\centering
 \includegraphics[scale=.4]{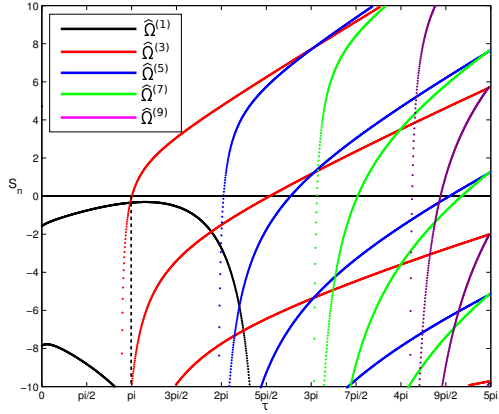}
 \caption{The $S_n$ map for different relative equilibria  with rotation frequency $\widehat{\Omega}^{(j)}(\tau)$. The zeroes of the $S_n$ map
  determine imaginary roots of $P_U(\lambda,\tau)$ in \eqref{eq:Nn_charact_functions}.}
 \label{fig:Sn_map_Pj}
\end{figure}

For this analysis we shall follow the steps described in section~\ref{subsubsec:Fix_SN_Omega_neq_0} for bifurcations in $\Fix(\S_N)$ and we will continue with the same example. We look for bifurcation points of the characteristic function $P_U$ with $\mu=1$, $K=1$, $\omega_M=1$ choosing $N=2$. Using the $S_n$ map as before, bifurcation points are shown in figure~\ref{fig:Sn_map_Pj}.
When $S_n=0$ then relative Hopf bifurcation of  non-synchronized RPOs occurs. In all cases $\delta(\omega(\tau^*))\geq 0$, i.e., roots cross the imaginary axis from the left to the right. The first bifurcation appears at $\tau=\pi$, which is lower than the lowest bifurcation value in $\Fix(\S_N)$ in figure~\ref{fig:Sn_map_Pfix}.


\section{Phase-difference model}
\label{sec:phase_diff_model}
In~\cite{Bueno2010, Correa2013, Bueno2009} an alternative approach is used to model a fully connected PLL network
using the phase difference between any two nodes $j$ and $k$ which is defined as 
\begin{equation}\label{e.defPhaseDifferences}
\varphi^{(j,k)}(t):=\theta_j(t)-\theta_{k}(t-\tau).
\end{equation}
From \eqref{eq:N-node_2nd-order_network} we get
\begin{equation}
  \label{eq:N-node_second-order_network}
  \ddot{\varphi}^{(i,j)}+\mu\dot{\varphi}^{(i,j)}+\frac{K\mu}{N-1}\left[\sum_{\stackrel[l\neq i]{l=1}{}}^N\sin\left(\varphi^{(i,l)}+\omega_M\tau \right)-\sum_{\stackrel[l\neq j]{l=1}{}}^N\sin\left(\varphi^{(j,l)}_{\tau}+\omega_M\tau\right) \right]=0.
\end{equation}
This phase-difference model has $\S_N$-symmetry as it is derived from the phase model of section~\ref{sec:Phase_model}, but not translational symmetry
 due to the definition of phase differences.

If $N>2$ then the phase difference model lives in $\varphi \in \R^{N(N-1)}$ and
$N(N-1)>  N$ so that it is to be expected that the phase difference model has fictitious solutions which do not
correspond to the phase model for $\theta \in \R^N$ from Section \ref{sec:Phase_model}.

 If $N=2$ and  $\theta_1([-\tau,0))$ and $\theta_2([-\tau,0))$ are known we can compute $\varphi^{(1,2)}$ and $\varphi^{(2,1)}$ for $t=[-\tau,\infty)$, and then compute $\theta(t)$ for $\forall t$ using~\eqref{e.defPhaseDifferences}. But from solutions $(\varphi^{(1,2)},\varphi^{(2,1)})$ we can not reconstruct $(\theta_1,\theta_2)$ initial data, so any solution of the $(\theta_1,\theta_2)$-dynamics is a solution of $(\varphi^{(1,2)},~\varphi^{(2,1)})$-dynamics, but not vice versa. Hence not all solutions of the $(\varphi^{(1,2)},~\varphi^{(2,1)})$-dynamics give solutions of the original $(\theta_1,\theta_2)$-dynamics.
We demonstrate this issue by studying equilibria of the phase-difference model.


In $\Fix(\S_N)$  an $N$-node network modelled by \eqref{eq:N-node_second-order_network}  has
a family of $\S_N$-invariant equilibria given by 
\begin{equation}
\varphi^{(i,j)}(t)=C,~~~C\in\R,
\label{eq:2n_phase_difference_const}
\end{equation}
which implies 
\begin{equation}
\label{eq:relative_periodic_sol}
\theta_i(t)-\theta_j(t-\tau)=\theta_j(t)-\theta_i(t-\tau)=C~~
\mbox{for all}~~
i\neq j. 
\end{equation}
The solutions are relative periodic orbits 
(RPOs) of the $\theta$-dynamics. If $N=2$ then 
\begin{equation}
\theta_i(t)=\theta_i(t-2\tau)+2C, i=1,2, 
\end{equation}
and if $C\equiv 0$ then $\theta(t)$ is $2\tau$-periodic with $\Z_2$-spatio-temporal symmetry, otherwise
the RPOs have $\Z_2$-spatio-temporal symmetry in a suitable comoving frame.

The $\theta$-dynamics from \eqref{eq:N-node_2nd-order_network} restricted to $\Fix(\S_N)$ is given by
\begin{equation}
  \label{eq:phaase_model_fixed_point_equ}
  \ddot{\theta}_i(t)+\mu\dot{\theta}_i(t)-K\mu\sin(\theta_i(t-\tau)-\theta_i(t)-\omega_M\tau)=0,\quad i=1,\ldots, N.
\end{equation}
By substituting \eqref{eq:2n_phase_difference_const} into   \eqref{eq:phaase_model_fixed_point_equ} we obtain the second-order ODE
\begin{equation}
  \ddot{\theta}_i(t)+\mu\dot{\theta}_i(t)+K\mu\sin(C+\omega_M\tau)=0,\quad i=1,\ldots, N
\end{equation}
whose solution is
\begin{equation}
  \theta_i(t)=-K\sin(C+\omega_M\tau) t+  C_1 +C_2e^{-\mu t},\quad i=1,\ldots, N,
\end{equation}
with  arbitrary constants $C,C_1,C_2$.
\marginnote{
\tiny
\begin{align*}
x &=\dot\theta_i, ~~a= K \mu \sin(C+\omega_M \tau).\\
  \dot x &= - \mu x -a\\
x(t) &= \e^{-t \mu} x_0 -a \int_0^t \e^{-\mu(t-s)} ds\\
&  = \e^{-t \mu} x_0 -\frac{a \e^{-\mu t}}{\mu} ( \e^{t \mu}-1 ) \\
&=\e^{-t \mu} x_0 -\frac{a }{\mu}  +\frac{a\e^{t \mu} }{\mu}\\
\theta_i(t) &= -\e^{-t \mu}\frac{ x_0}{\mu} - \frac{a t }{\mu}  + a\e^{t \mu} + const
\end{align*}
\normalsize
}
 From \eqref{eq:relative_periodic_sol} we have for any $i\neq j$
\begin{equation}
C=  \theta_i(t)-\theta_j(t-\tau)=-K\tau\sin(C+\omega_M\tau)+C_2e^{-\mu t}(1-e^{\mu\tau}).
\end{equation}
Hence, $C_2=0$ and so 
\begin{equation}
  C=-K\tau\sin(C+\omega_M\tau),
\end{equation}
and so $C = \Omega(\tau) \tau$ with 
  $\Omega(\tau)$  from \eqref{eq:Omega(tau)_definition}.

So although the N-node model from ~\eqref{eq:N-node_second-order_network} admits an $\S_N$-invariant 
equilibrium $\varphi_{ij}(t)\equiv C$  for all $C\in \R$ only the choice $C=\Omega(\tau)\tau$ corresponds to an actual $\S_N$-invariant 
equilibrium of the phase model from  section~\ref{sec:Phase_model}.  

The matrix $L$ from \eqref{eq:L(tau)_definition}
of  the linearization around
 the equilibrium given by $\varphi^{(ij)}(t)\equiv C$ for all $i\neq j$ of   \eqref{eq:N-node_second-order_network} for $N=2$ is given by 
\begin{equation}
  L=\left(
    \begin{array}{cccc}
      0&1&0&0\\
      -K\mu\cos(C+\omega_M\tau)&-\mu&K\mu\cos(C+\omega_M\tau)e^{-\lambda\tau}&0\\
      0&0&0&1\\
      K\mu\cos(C+\omega_M\tau)e^{-\lambda\tau}&0&-K\mu\cos(C+\omega_M\tau)&-\mu
    \end{array}
\right).
\end{equation}
The characteristic matrix $\triangle(\lambda,\tau):=\lambda I -L$ can be uncoupled into blocks corresponding to
isotypic components giving
\begin{equation}
  \rho\triangle\rho^{-1}=\left(\begin{array}{c|c}\triangle_1&0\\\hline 0&\triangle_2\end{array}\right)
\end{equation}
for some transformation $\rho \in \Mat(4)$, and the characteristic functions 
$P_j(\lambda,\tau)=\det(\triangle_j(\lambda,\tau))$, $j=1,2$, are
\begin{equation}
  \begin{array}{l}
    P_{1,2}(\lambda,\tau)=\lambda^2+\mu\lambda\mp K\mu\cos(C+\omega_M\tau)+ K\mu\cos(C+\omega_M\tau)e^{-\lambda\tau},
  \end{array}
\end{equation}
which upon substituting $C=\Omega(\tau)\tau$, gives the characteristic functions  for the phase model in \eqref{eq:Nn_charact_functions} with $N=2$, where $\Omega(\tau)\tau+\omega_M\tau=\widehat{\Omega}\tau$.

For a 3-node network using the phase differences model linearizing around the equilibrium point $\varphi^{(i,j)}=\Omega(\tau)\tau$ we obtain
\begin{equation}
\begin{array}{r}
\det \triangle(\lambda,\tau)=(\lambda^2+\mu\lambda)^3(\lambda^2+\mu\lambda+K\mu\cos(\widehat{\Omega}\tau)-K\mu\cos(\widehat{\Omega}\tau)e^{-\lambda\tau})\\\\
(\lambda^2+\mu\lambda+K\mu\cos(\widehat{\Omega}\tau)+\dfrac{K\mu}{2}\cos(\widehat{\Omega}\tau)e^{-\lambda\tau})^2=0,
\end{array}
\end{equation}
\note{DF: We Know that for
  equation~\ref{eq:N-node_second-order_network}, the linearization at
  equilibria $\varphi=C\in\R$, restricted to the $\varphi^{(i,j)}$
  component is:
  \begin{align*}
    (L)(x)|_{x^{(i,j)}}=\left(
      \begin{array}{c}
        x_2^{(i,j)}\\
        -\mu x_2^{(i,j)}-a\sum_{\stackrel[l\neq
          i]{l=1}{}}^Nx_1^{(i,l)}+a\e^{-\lambda\tau}\sum_{\stackrel[l\neq
          j]{l=1}{}}^Nx_1^{(j,l)}
      \end{array}
\right);
\end{align*}
$x_1^{(i,j)}=\varphi^{(i,j)}$, $x_2^{(i,j)}=\dot\varphi^{(i,j)}$,
$x^{(i,j)}=[x_1^{(i,j)},x_2^{(i,j)}]^T$,
$a=\dfrac{K\mu}{N-1}\cos(C+\omega_M\tau)$.
For $N=3$, $L$ has the form:
\begin{align*}
  L=\left(
  \begin{array}{ccc}
    L_x&L_1&L_2\\
    L_2&L_x&L_1\\
    L_1&L_2&L_x
  \end{array}\right),
\end{align*}
\begin{align*}
  L_x=\left(
  \begin{array}{cccc}
    0&1&0&0\\
    -a&\mu&-a&0\\
    0&0&0&1\\
    -a&0&-a&-\mu
  \end{array}\right),\\L_1=\left(
  \begin{array}{cccc}
    0&0&0&0\\
    a\e^{-\lambda\tau}&0&a\e^{-\lambda\tau}&0\\
    0&0&0&0\\
    0&0&0&0
  \end{array}\right),\\L_2=\left(
  \begin{array}{cccc}
    0&0&0&0\\
    0&0&0&0\\
    0&0&0&0\\
    a\e^{-\lambda\tau}&0&a\e^{-\lambda\tau}&0
  \end{array}\right).
\end{align*}
Since this a determinant computation, I used Matlab code:
syms K u tau C  a real\\
syms La\\
Lx = [ 0  1  0  0;
      -a -u -a  0;
       0  0  0  1;
      -a  0 -a -u];\\
L1 = [0              0             0  0;
      a*exp(-La*tau) 0 a*exp(-La*tau) 0;
      0              0             0  0;
      0              0             0  0];\\
L2 = [0              0             0  0;
      0              0             0  0;
      0              0             0  0;
      a*exp(-La*tau) 0 a*exp(-La*tau) 0];\\  
L = [Lx L1 L2;
     L2 Lx L1;
     L1 L2 Lx];\\  
M = eye(12)*La - L;\\  
detM = simplify(det(M));
\begin{align*}
\begin{array}{l}
 \text{det}(\text{I}\lambda-L)=\\{\lambda}^3\, \e^{- 3\, \lambda\,
   \tau}\, {\left(\lambda + u\right)}^3\, {\left(a + 2\, a\,
     \e^{\lambda\, \tau} + {\lambda}^2\, \e^{\lambda\, \tau} +
     \lambda\, u\, \e^{\lambda\, \tau}\right)}^2\, \left(2\, a\,
   \e^{\lambda\, \tau} - 2\, a + {\lambda}^2\, \e^{\lambda\, \tau} +
   \lambda\, u\, \e^{\lambda\, \tau}\right)=0
\end{array}
\end{align*}
}
\note{Rearranging terms we get:
  \begin{align*} \text{det}(\text{I}\lambda-L)=(\lambda^2+\mu\lambda)^3\left(\lambda^2+\lambda\mu+2a+a\e^{-\lambda\tau}\right)^2\left(\lambda^2+\lambda\mu+2a-2a\e^{-\lambda\tau}\right)=0
  \end{align*}
Substituting $a=\dfrac{K\mu}{2}\cos(C+\omega_M\tau)$ we finally get:
  \begin{align*}
\begin{array}{l}
  \text{det}(\text{I}\lambda-L)=\\(\lambda^2+\mu\lambda)^3\left(\lambda^2+\lambda\mu+K\mu\cos(C+\omega_M\tau)+\dfrac{K\mu}{2}\cos(C+\omega_M\tau)\e^{-\lambda\tau}\right)^2\\\left(\lambda^2+\lambda\mu+K\mu\cos(C+\omega_M\tau)-K\mu\cos(C+\omega_M\tau)\e^{-\lambda\tau}\right)=0.
\end{array}
  \end{align*}
}
and here is clear that the term $(\lambda^2+\mu\lambda)^3$ does not correspond to roots of the characteristic function  for the phase model from \eqref{eq:Nn_charact_functions} with $N=3$.


 
\section{Discussion and conclusions}
\label{sec:conclusions}
Due to the $\S_N$-symmetry of a second-order N-node oscillators network modeled using the full-phase variables
the linearization along every $\S_N$-invariant equilibrium has $N$ blocks; one  of them corresponds 
to the fixed point space $\Fix(\S_N)$ and the others are identical and lead to eigenvalues of multiplicity $N-1$. This decomposition simplifies the bifurcation analysis considerably. Hopf bifurcation in the first block is symmetry-preserving, i.e., leads to bifurcation
of synchronized periodic orbits which are periodic orbits in  $\Fix(\S_N)$. Hopf bifurcation in the other blocks is symmetry-breaking and leads to bifurcation of non-synchronized or partially synchronized periodic orbits.

We presented  this decomposition  for  second order oscillators,  but networks of  higher order oscillators, modelled accordingly,  also have
$\S_N$-symmetry, therefore the above mentioned block decomposition
also applies provided the time-delay between the nodes are the same.

Although the full-phase model is obtained from the phase model by including the double-frequency term
and passing into a comoving frame with velocity $\omega_M$, the
dynamics observed in each case are significantly different
due to the translation symmetry
of the phase model which is not present in the full-phase model.
The translation symmetry causes the existence of relative equilibria
or circles of equilibria if the rotation frequency vanishes. As a consequence, the   linearization of the relative equilibria in the corotating 
frame always has a zero eigenvalue. At the relative Hopf bifurcations discussed in Section \ref{sec:Phase_model} relative periodic
orbits emanate which correspond to quasiperiodic orbits in the original frame
of reference. Due to the translation symmetry frequency locking does not appear
on those invariant tori.

The phase-difference model discussed in section~\ref{sec:phase_diff_model} introduces fictitious solutions that may not correspond to real solutions of the phase model analyzed in section~\ref{sec:Phase_model} even when the equilibrium point is chosen as $\Omega(\tau)\tau$ for $N>2$. 
 
We conclude that of the three models studied here, only the full-phase model from section~\ref{sec:Full_phase_model} represents better and without any approximations the dynamics of a fully connected N-node time-delay network.

The  stability  of  the bifurcating periodic solutions  will be the focus of future work.

\section*{Acknowledgement}
We would like to thank   UGPN,   the Department of Mathematics of the University of Surrey and  Escola Polit\'ecnica da Universidade de S\~ao Paulo for their support.

\end{document}